\documentclass[12pt,draftcls,onecolumn,twoside]{IEEEtran}
\usepackage[latin1]{inputenc}
\usepackage[mathcal]{euscript}
\usepackage{amssymb,amsmath,amsthm}
\usepackage{graphicx}
\usepackage[english]{babel}
\newcommand{\Pmp}{\mathbb{P}_{\check{\mathrm{M}}\mathrm{P}}}

\newtheorem{remark}{Remark}
\newtheorem{prop}{Proposition}

\newtheorem{theo}{Theorem}
\newtheorem{lemma}{Lemma}

\newtheorem{assumption}{Assumption}

\usepackage{color}
\definecolor{MyDarkBlue}{rgb}{0.1,0,0.65}
\definecolor{MyDarkRed}{rgb}{.85,0,0.1}

\newcommand{\tr}{\mathrm{tr}}

\newcommand{\f}{\mathbf{f}}
\newcommand{\F}{\mathbf{F}}
\newcommand{\tildef}{\tilde{\mathbf{f}}}
\newcommand{\lspiked}{\lambda_{\mathrm{spk}}^\infty}

\newcommand{\RR}{\vspace*{.2cm}\\}

\newcommand{\ud}{\textrm{d}}

\newcommand{\bs}{\boldsymbol}
\newcommand{\dsp}{\displaystyle}

\newcommand{\bb}{\mathbb}

\newcommand{\gtrlessHO}{
  \begin{array}{c}
    {\tiny{H_0}}\\ {\gtrless} \\{\tiny H_1}
  \end{array}}
\newcommand{\gtrlessHU}{
  \begin{array}{c}
    {\footnotesize H_1}\\ {\gtrless} \\{\footnotesize H_0}
  \end{array}}

\title{Performance of Statistical Tests for Single Source Detection using Random Matrix Theory}

\author{P. Bianchi, M. Debbah, M. Maida and J. Najim\\
  \today \thanks{This work was partially supported by french programs
    ANR-07-MDCO-012-01 `Sesame', and ANR-08-BLAN-0311-03
    `GranMa'.}
  \thanks{P. Bianchi and J. Najim are with CNRS and Télécom Paristech,
    France. {\tt \{bianchi,najim\}@telecom-paristech.fr}\ ,}
  \thanks{M. Debbah is with SUPELEC and holds Alcatel-Lucent/Supélec
    Flexible Radio chair, France {\tt merouane.debbah@supelec.fr}\ ,}
  \thanks{M. Maida is with Université Paris-Sud, UMR CNRS 8628,
    France. {\tt Mylene.Maida@math.u-psud.fr}\ ,}%
}

\begin{document}
\maketitle

\begin{abstract}
  This paper introduces a unified framework for the detection of a single
  source with a sensor array in the context where
  the noise variance and the channel between the source and the
  sensors are unknown at the receiver. The Generalized Maximum
  Likelihood Test is studied and yields the analysis of the ratio
  between the maximum eigenvalue of the sampled covariance matrix and
  its normalized trace. Using recent results from random matrix theory,
  a practical way to evaluate the threshold and the $p$-value of the
  test is provided in the asymptotic regime where the number $K$ of
  sensors and the number $N$ of observations per sensor are large but
  have the same order of magnitude. The theoretical performance of the
  test is then analyzed in terms of Receiver Operating Characteristic
  (ROC) curve. It is in particular proved that both Type~I and Type~II
  error probabilities converge to zero exponentially as the dimensions
  increase at the same rate, and closed-form expressions are provided
  for the error exponents. These theoretical results rely on a precise
  description of the large deviations of the largest eigenvalue of
  spiked random matrix models, and establish that the presented test
  asymptotically outperforms the popular test based on the condition
  number of the sampled covariance matrix.

\end{abstract}

\section{Introduction}
The detection of a source by a sensor array is at the heart of many
wireless applications. It is of particular interest in the realm of
cognitive radio \cite{MiMa99,Ha05} where a multi-sensor cognitive
device (or a collaborative network\footnote{The collaborative network
  corresponds to multiple base stations connected, in a wireless or
  wired manner, to form a virtual antenna system\cite{DoLeAg02}.})
needs to discover or sense by itself the surrounding environment. This
allows the cognitive device to make relevant choices in terms of
information to feed back, bandwidth to occupy or transmission power to
use.
When the cognitive device is switched on, its prior knowledge (on the
noise variance for example) is very limited and can rarely be
estimated prior to the reception of data. This unfortunately rules out
classical techniques based on energy detection
\cite{Ur67,Ko02,SiDiAl03} and requires new sophisticated techniques
exploiting the space or spectrum dimension.

In our setting, the aim of the multi-sensor cognitive detection phase
is to construct and analyze tests associated with the following
hypothesis testing problem:
\begin{equation}\label{eq_hypothesis_test}
{\bs y}(n)= \left\lbrace\begin{array}{ll}
{\bs w}(n)  & \textrm{under}\ H_0\RR
\dsp {\bs h}\:{s}(n) + {\bs w}(n)  & \textrm{under}\ H_1\\
\end{array}
\right. \qquad \textrm{for}\quad n=0:N-1\ ,
\end{equation}
where ${\bs y}(n) = [y_1(n),\dots,y_K(n)]^T$ is the observed $K× 1$
complex time series
, ${\bs w}(n)$ represents a
$K× 1$ complex circular Gaussian white noise process with unknown variance 
$\sigma^2$, and $N$ represents the number of received samples. Vector ${\bs
  h}\in{\mathbb{C}}^{K× 1}$ is a deterministic vector and typically
represents the propagation channel between the source and the $K$
sensors. Signal $s(n)$ denotes a standard scalar independent and identically distributed (i.i.d.) circular complex
Gaussian process with respect to the samples $n=0:N-1$ and stands for
the source signal to be detected.

The standard case where the propagation channel and the noise variance
are known has been thoroughly studied in the literature in the Single
Input Single Output case \cite{Ur67,Ko02,SiDiAl03} and Multi-Input
Multi-Ouput \cite{QuCuSaPo08} case.  In this simple context, the most
natural approach to detect the presence of source $s(n)$ is the
well-known \emph{Neyman-Pearson} (NP) procedure which consists in
rejecting the null hypothesis when the observed likelihood ratio lies
above a certain threshold \cite{lehman}.  Traditionally, the value of
the threshold is set in such a
way 
that the \emph{Probability of False Alarm} (PFA) is no larger than a
predefined \emph{level} $\alpha\in(0,1)$.  Recall that the PFA
(resp. the miss probability) of a test is defined as the probability
that the receiver decides hypothesis $H_1$ (resp. $H_0$) when the true
hypothesis is $H_0$ (resp.  $H_1$). The NP test is known to be
uniformly most powerful \emph{i.e.}, for any level $\alpha\in(0,1)$,
the NP test has the minimum achievable miss probability (or
equivalently the maximum achievable power) among all tests of level
$\alpha$. In this paper, we assume on the opposite that:
\begin{itemize}
\item the noise variance $\sigma^2$ is unknown,
\item vector $\bs h$ is unknown.
\end{itemize}
In this context, probability density functions of the observations
${\bs y}(n)$ under both $H_0$ and $H_1$ are unknown, and the classical
NP approach can no longer be employed.  As a consequence, the
construction of relevant tests for~(\ref{eq_hypothesis_test}) together
with the analysis fo their perfomances is a crucial issue.  The
classical approach followed in this paper consists in replacing the
unknown parameters by their maximum likelihood estimates. This leads
to the so-called \emph{Generalized Likelihood Ratio} (GLR). The
\emph{Generalized Likelihood Ratio Test} (GLRT), which rejects the
null hypothesis for large values of the GLR, easily reduces to the
statistics given by the ratio of the largest eigenvalue of the sampled
covariance matrix with its normalized trace,
cf. \cite{wax85,ChBoMeHe07,Gaz10}. Nearby
statistics~\cite{MesIT08,MesSP08,KRNA08,KRNA09}, with good practical
properties, have also been developed, but would not yield a different
(asymptotic) error exponent analysis.

In this paper, we analyze the performance of the GLRT in the
asymptotic regime where the number $K$ of sensors and the number $N$
of observations per sensor are large but have the same order of
magnitude. This assumption is relevant in many applications, among
which cognitive radio for instance, and casts the problem into a large
random matrix framework.

Large random matrix theory has already been applied to signal
detection \cite{SilCom92} (see also \cite{CouDeb10pre}), and recently to hypothesis testing
\cite{KRNA09,RaoSpe08,RaoSil09pre}. In this article, the focus is
mainly devoted to the study of the largest eigenvalue of the sampled
covariance matrix, whose behaviour changes under $H_0$ or $H_1$. The
fluctuations of the largest eigenvalue under $H_0$ have been described
by Johnstone \cite{Joh01} by means of the celebrated Tracy-Widom
distribution, and are used to study the threshold and the $p$-value of
the GLRT.

In order to characterize the performance of the test, a natural
approach would have been to evaluate the \emph{Receiver Operating
  Characteristic} (ROC) curve of the GLRT, that is to plot the power
of the test versus a given level of confidence. Unfortunately, the ROC
curve does not admit any simple closed-form expression for a finite
number of sensors and snapshots. As the miss probability of the GLRT
goes exponentially fast to zero, the performance of the GLRT is
analyzed via the computation of its error exponent, which caracterizes
the speed of decrease to zero. Its computation relies on the study of
the large deviations of the largest eigenvalue of 'spiked' sampled
covariance matrix. By 'spiked' we refer to the case where the
eigenvalue converges outside the bulk of the limiting spectral
distribution, which precisely happens under hypothesis $H_1$. We build
upon \cite{Mai07} to establish the large deviation principle, and
provide a closed-form expression for the rate function.  

We also introduce the error exponent curve, and plot the error
exponent of the power of the test versus the error exponent for a
given level of confidence. The error exponent curve can be interpreted
as an asymptotic version of the ROC curve in a $\log$-$\log$ scale and
enables us to establish that the GLRT outperforms another test based
on the condition number, and proposed
by~\cite{zengliang,art-cardoso-2008,garello} in the context of
cognitive radio.

Notice that the results provided here (determination of the threshold
of the GLRT test and the computation of the error exponents) would
still hold within the setting of real Gaussian random variables
instead of complex ones, with minor modifications\footnote{Details are provided in Remarks \ref{rem:real1}
and \ref{rem:real2}.}.

The paper is organized as follows. 

Section~\ref{sec:glrt} introduces the GLRT. The value of the
threshold, which completes the definition of the GLRT, is established
in Section \ref{sec:exact}. As the latter threshold has no simple
closed-form expression and as its practical evaluation is difficult,
we introduce in Section~\ref{subsec:asymptotic} an asymptotic
framework where it is assumed that both the number of sensors $K$ and
the number $N$ of available snapshots go to infinity at the same
rate. This assumption is valid for instance in cognitive radio
contexts and yields a very simple evaluation of the threshold, which
is important in real-time applications.

In Section \ref{sec:known}, we recall several results of large random
matrix theory, among which the asymptotic fluctuations of the largest
eigenvalue of a sample covariance matrix, and the limit of the largest
eigenvalue of a spiked model.
 
These results are used in Section \ref{sec:asympt} where an
approximate threshold value is derived, which leads to the same PFA as
the optimal one in the asymptotic regime. This analysis yields a
relevant practical method to approximate the $p$-\emph{values}
associated with the GLRT.

Section~\ref{sec:power} is devoted to the performance analysis of the
GLRT. We compute the error exponent of the GLRT, derive its expression
in closed-form by establishing a \emph{Large Deviation Principle}
for the test statistic $T_N$ \footnote{Note that in recent
  papers \cite{ABMASA09,KRNA08,KRNA09}, the fluctuations of the test
  statistics under $H_1$, based on large random matrix techniques,
  have also been used to approximate the power of the test. We believe
  that the performance analysis based on the error exponent approach,
  although more involved, has a wider range of validity.}, and
describe the error exponent curve.

Section~\ref{sec:condition} introduces the test based on the condition
number, that is the statistics given by the ratio between the largest
eigenvalue and the smallest eigenvalue of the sampled covariance
matrix.  We provide the error exponent curve associated with this test
and prove that the latter is outperformed by the GLRT.

Section~\ref{sec:numerical} provides further numerical illustrations and
conclusions are drawn in Section~\ref{sec:conclusion}.

Mathematical details are provided in the Appendix. In particular, a
full rigorous proof of a large deviation principle is provided in
Appendix \ref{app:ldp}, while a more informal proof of a nearby large
deviation principle, maybe more accessible to the non-specialist, is
provided in Appendix \ref{app:ldp-U}.

\subsection*{Notations}

For $i\in\{0,1\}$, ${\mathbb P}_{i}({\cal E})$ represents the
probability of a given event $\cal E$ under hypothesis $H_i$.  For any
real random variable $T$ and any real number $\gamma$, notation
$$
T\:
_{_{H_0}}\!\!\gtrless^{^{H_1}} \gamma
$$
stands for the test function which rejects the null hypothesis when
$T>\gamma$. In this case, the \emph{probability of false alarm (PFA)}
of the test is given by ${\bb P}_0(T>\gamma)$, while the power of the
test is ${\bb P}_1(T>\gamma)$. Notation $\xrightarrow[H_i]{a.s.}$
stands for the almost sure (a.s.) convergence under hypothesis
$H_i$. For any one-to-one mapping $T:{\cal X}\to {\cal Y}$ where $\cal
X$ and $\cal Y$ are two sets, we denote by $T^{-1}$ the inverse of $T$
w.r.t. composition. For any borel set $A\in{\bb R}$, $x\mapsto {\bs
  1}_A(x)$ denotes the indicator function of set $A$ and $\|{\bs x}\|$
denotes the Euclidian norm of a given vector $\bs x$.  If ${\bs A}$ is
a given matrix, denote by ${\bs A}^H$ its transpose-conjugate.  If $F$
is a cumulative distribution function (c.d.f.), we denote by $\bar{F}$
is complementary c.d.f., that is: $\bar{F}= 1-F$.


\section{Generalized Likelihood Ratio Test}
\label{sec:glrt}
In this section, we derive the Generalized Likelihood Ratio Test
(section \ref{subsec:glrt}) and compute the associated threshold and
$p$-value (section \ref{sec:exact}).  This exact computation raises
some computational issues, which are circumvented by the introduction
of a relevant asymptotic framework, well-suited for mathematical
analysis (Section \ref{subsec:asymptotic}).

\subsection{Derivation of the Test}\label{subsec:glrt}
Denote by $N$ the number of observed samples and recall that:
$$
{\bs y}(n)= \left\lbrace\begin{array}{ll}
{\bs w}(n)  & \textrm{under}\ H_0\RR
\dsp {\bs h}\:{s}(n) + {\bs w}(n)  & \textrm{under}\ H_1\\
\end{array}
\right.,\qquad n=0:N-1\ ,
$$
where $({\bs w}(n),0\leq n\leq N-1)$ represents an independent and
identically distributed (i.i.d.) process of $K× 1$ vectors with
circular complex Gaussian entries with mean zero and covariance matrix
$\sigma^2 {\bf I}_K$, vector ${\bs h}\in{\mathbb{C}}^{K× 1}$ is
deterministic, signal $(s(n),0\leq n\leq N-1)$ denotes a scalar i.i.d.
circular complex Gaussian process with zero mean and unit variance.
Moreover, $({\bs w}(n),0\leq n\leq N-1)$ and $(s(n),0\leq n\leq N-1)$ are
assumed to be independent processes. We stack the observed data into a
$K× N$ matrix ${\bf Y} = [{\bs y}(0),\dots,{\bs y}(N-1)]$. Denote by
$\hat{\bf R}$ the sampled covariance matrix:
$$
\hat {\bf R} = \frac{1}{N} {\bf Y}{\bf Y}^H,
$$
and respectively, by $p_{0}({\bf Y} ; \sigma^2)$ and $p_{1}({\bf
  Y} ; {\bs h}, \sigma^2)$ the likelihood functions of the observation
matrix ${\bf Y}$ indexed by the unknown parameters ${\bs h}$ and
$\sigma^2$ under hypotheses $H_0$ and $H_1$.

As ${\bf Y}$ is a $K × N$ matrix whose columns are i.i.d.
Gaussian vectors with covariance matrix ${\bf \Sigma}$ defined by:
\begin{equation}\label{eq:covariance}
{\bf \Sigma} =\left\{
\begin{array}{ll}
\sigma^2 {\bf  I}_K&\textrm{under}\ H_0\\
{\bs h}{\bs h}^H + \sigma^2 {\bf I}_K& \textrm{under}\ H_1
\end{array}\right. \ ,
\end{equation}
the likelihood functions write:
\begin{align}
\label{eq:pzero}
  p_{0}({\bf Y} ; \sigma^2) &= (\pi \sigma^2)^{-NK} \exp\left( -\frac{N}{\sigma^2}\tr \: {\hat{\bf R}} \right)\ , \\
  p_{1}({\bf Y} ; {\bs h},\sigma^2 ) &= (\pi^K \det({\bs h}{\bs h}^H + \sigma^2 {\bf I}_K))^{-N} \exp\left(
    -N\tr\:(\hat {\bf R} ({\bs h}{\bs h}^H + \sigma^2 {\bf I}_K)^{-1})\right)\ .
\label{eq:pun}
\end{align}

In the case where parameters $\bs h$ and $\sigma^2$ are available, the
celebrated Neyman-Pearson procedure yields a uniformly most powerful
test, given by the likelihood ratio statistics $\frac{p_{1}({\bf Y} ;
  {\bs h}, \sigma^2)}{p_{0}({\bf Y} ; \sigma^2 )}$.

However, in the case where $\bs h$ and $\sigma^2$ are unknown, which
is the problem addressed here, no simple procedure
garantees a uniformly most powerful test, and a classical approach
consists in computing the GLR:
\begin{equation}
  \label{eq:LR}
  L_N  = \frac{\sup_{{\bs h}, \sigma^2} p_{1}({\bf Y} ; {\bs h},\sigma^2 )}{\sup_{\sigma^2} p_{0}({\bf Y} ; \sigma^2)}\:.
\end{equation}
In the GLRT procedure, one rejects
hypothesis $H_0$ whenever $L_N>\xi_N$, where $\xi_N$ is a certain
threshold which is selected in order that the PFA
${\bb P}_0(L_N>\xi_N)$ does not exceed a given level $\alpha$.

In the following proposition, which follows after straightforward
computations from~\cite{and63} and \cite{wax85}, we derive the closed form
expression of the GLR $L_N$. Denote by $\lambda_1 > \lambda_2 >\dots
>\lambda_K\geq0$ the ordered eigenvalues of $\hat{\bf R}$ (all
distincts with probability one).
\begin{prop}
\label{prop:glr}
Let $T_N$ be defined by:
\begin{equation}
  \label{eq:mu}
  T_N = \frac{\lambda_1}{\frac{1}{K}\:\tr\:\hat{\bf R}}\:,
\end{equation}
then, the GLR (cf. Eq.~(\ref{eq:LR})) writes:
$$
L_N = \frac{C}{\left(T_N\right)^N\left(1-\frac{T_N}{K}\right)^{(K-1)N}}
$$
where $C=\left(1-\frac{1}{K}\right)^{(1-K)N}$.
\end{prop}
By Proposition~\ref{prop:glr}, $L_N=\phi_{N,K}(T_N)$ where $\phi_{N,K} : x\mapsto
Cx^{-N}\left(1-\frac{x}{K}\right)^{N(1-K)}$. The GLRT rejects the null
hypothesis when inequality $L_N>\xi_N$ holds. As $T_N\in (1,K)$ with probability one and as
$\phi_{N,K}$ is increasing on this interval, the latter inequality is equivalent to
$T_N>\phi_{N,K}^{-1}(\xi_N)$. Otherwise stated, the GLRT reduces to
the test which rejects the null hypothesis for large values of~$T_N$:
\begin{equation}
  \label{eq:TN}
  T_N \gtrlessHU \gamma_N
\end{equation}
where $\gamma_N=\phi_{N,K}^{-1}(\xi_N)$ is a certain threshold which
is such that the PFA does not exceed a given
level $\alpha$. In the sequel, we will therefore focus on the
test statistics $T_N$.

\begin{remark} There exist several variants of the above
  statistics~\cite{MesIT08,MesSP08,KRNA08,KRNA09}, which merely
  consist in replacing the normalized trace with a more involved
  estimate of the noise variance. Although very important from a
  practical point of view, these variants have no impact on the
  (asymptotic) error exponent analysis. Therefore, we restrict our
  analysis to the traditional GLRT for the sake of simplicity.
\end{remark}

\subsection{Exact threshold and $p$-values}
\label{sec:exact}

In order to complete the construction of the test, we must provide a
procedure to set the threshold $\gamma_N$. As usual, we propose to
define $\gamma_N$ as the value which maximizes the power ${\bb
  P}_1(T_N>\gamma_N)$ of the test~(\ref{eq:TN}) while keeping the
PFA ${\bb P}_0(T_N>\gamma_N)$ under a desired
level $\alpha\in (0,1)$. It is well-known (see for instance
\cite{lehman,Van98}) that the latter threshold is obtained by:
\begin{equation}
  \gamma_N=p_N^{-1}(\alpha)\label{eq:optgamma}
\end{equation}
where $p_N(t)$ represents the complementary c.d.f. of the statistics
$T_N$ under the null hypothesis:
\begin{equation}
  \label{eq:ccdf}
  p_N(t) = {\bb P}_0(T_N>t)\:.
\end{equation}
Note that $p_N(t)$ is continuous and decreasing from 1 to 0 on
$t\in[0,\infty)$, so that the threshold $p_N^{-1}(\alpha)$ in~(\ref{eq:optgamma})
is always well defined.
When the threshold is fixed to $\gamma_N=p_N^{-1}(\alpha)$, the GLRT rejects
the null hypothesis when $T_N>p_N^{-1}(\alpha)$ or equivalently, when $p_N(T_N)<\alpha$.
It is usually convenient to rewrite the GLRT under the following form:
\begin{equation}
  \label{eq:testpval}
  p_N(T_N) \gtrlessHO \alpha\:.
\end{equation}
The statistics $p_N(T_N)$ represents the \emph{significance probability} or
$p$-value of the test.  The null hypothesis is rejected
when the $p$-value $p_N(T_N)$ is below the level $\alpha$.  In practice,
the computation of the $p$-value associated with one experiment is of
prime importance. Indeed, the $p$-value not only allows to accept/reject an
hypothesis by~(\ref{eq:testpval}), but it furthermore reflects how strongly the data
contradicts the null hypothesis~\cite{lehman}.

In order to evaluate $p$-values, we derive in the sequel the exact
expression of the complementary c.d.f. $p_N$. The crucial point is that $T_N$ is a function of the
eigenvalues $\lambda_1,\dots,\lambda_K$ of the sampled covariance matrix $\hat
{\bf R}$.
We have \begin{equation}
  \label{eq:ccdfexpr}
p_N(t) = \int_{\Delta_t}  p_{K,N}^{0}(x_1,\cdots,x_K)\ud x_{1:K}
\end{equation}
where for each $t$, the domain of integration $\Delta_t$ is defined by:
$$
\Delta_t=\left\{(x_1,\dots,x_K)\in \mathbb{R}^K,\ \frac{Kx_1}{x_1+\dots+x_K}>t\right\}\ ,
$$
and $p_{K,N}^{0}$ is the joint probability density function (p.d.f.) of
the ordered eigenvalues of ${\bf \hat R }$ under $H_0$ given by:
\begin{equation}
  \label{eq:exactpdf0}
  p_{K,N}^{0}(x_{1:K}) = \frac{{\bs 1}_{(x_1\geq\dots\geq x_K \geq 0)}}{ Z_{K,N}^{0}}\prod_{1\leq i<j\leq K}(x_j-x_i)^2 \:\prod_{j=1}^Kx_j^{N-K}
e^{-Nx_j}
\end{equation}
where ${\bs 1}_{(x_1\geq\dots\geq x_K \geq 0)}$ stands for the
indicator function of the set $\{(x_1\dots x_K)\: : \:
x_1\geq\dots\geq x_K \geq 0\}$ and where $ Z_{K,N}^{0}$ is the
normalization constant (see for instance \cite{Meh04}, \cite[Chapter
4]{AGZ09}).

\begin{remark}
  For each $t$, the computation of $p_N(t)$ requires the numerical
  evaluation of a non-trivial integral. 
  Despite the fact that powerful numerical methods, based on
  representations of such integrals with hypergeometric functions
  \cite{Mui82}, are available (see for instance \cite{Koev-toolbox},
  \cite{KoeEde06}), an \emph{on line} computation, requested in a
  number of real-time applications, may be out of reach.

  Instead, tables of the function $p_N$ should be computed \emph{off
    line} \emph{i.e.}, prior to the experiment. As both the dimensions
  $K$ and $N$ may be subject to frequent changes\footnote{In cognitive
    radio applications for instance, the number of users $K$ which are
    connected to the network is frequently varying.}, all possible
  tables of the function $p_N$ should be available at the detector's
  side, for all possible values of the couple $(N,K)$. This both
  requires substantial computations and considerable memory space. In
  what follows, we propose a way to overcome this issue.
\end{remark}

In the sequel, we study the asymptotic behaviour of the complementary
c.d.f. $p_N$ when both the number of sensors $K$ and the number of
snapshots $N$ go to infinity at the same rate. This analysis leads to
simpler testing procedure.

\subsection{Asymptotic framework}\label{subsec:asymptotic}

We propose to analyze the asymptotic behaviour of the complementary c.d.f. $p_N$ as
the number of observations goes to infinity. More precisely, we
consider the case where both the number $K$ of sensors and the number
$N$ of snapshots go to infinity at the same speed, as assumed below
\begin{equation}
  \label{eq:regime}
  N\to\infty,\; K\to\infty,\;\; c_N:= \frac KN\to c, \textrm{ with } 0<c<1.
\end{equation}
This asymptotic regime is relevant in cases where the sensing system
must be able to perform source detection in a moderate amount of time
\emph{i.e.}, the number $K$ of sensors and the number $N$ of samples
being of the same order. This is in particular the case in cognitive
radio applications (see for instance~\cite{phd-mitola-2000}). Very
often, the number of sensors is lower than the number of snapshots,
hence the ratio $c$ lower than 1.

In the sequel, we will simply denote $N,K\to \infty$ to refer
to the asymptotic regime \eqref{eq:regime}.

\begin{remark} The results related to the GLRT presented in Sections
  \ref{sec:asympt} and \ref{sec:power} remain true for $c\ge 1$; in
  the case of the test based on the condition number and presented in
  Section \ref{sec:condition}, extra-work is needed to handle the fact
  that the lowest eigenvalue converges to zero, which happens if $c\ge
  1$.
\end{remark}

\section{Large random matrices - Largest eigenvalue - Behaviour of the GLR statistics}
\label{sec:known}

In this section, we recall a few facts on large random matrices as the
dimensions $N,K$ go to infinity. We focus on the behaviour of the
eigenvalues of $\bf \hat R$ which differs whether hypothesis $H_0$
holds (Section \ref{subsec:H0}) or $H_1$ holds (Section
\ref{subsec:H1}).

As the column vectors of ${\bf Y}$ are i.i.d. complex Gaussian with covariance matrix
${\bf \Sigma}$ given by \eqref{eq:covariance}, the probability density of ${\bf \hat R}$
is given by:
$$
\frac{1}{Z(N,K,{\bf \Sigma})} e^{-N \tr ({\bf \Sigma}^{-1}{\bf \hat R} )} (\textrm {det } {\bf \hat R})^{N-K},
$$
where $Z(N,K,{\bf \Sigma})$ is a normalizing constant.

\subsection{Behaviour under hypothesis $H_0$}\label{subsec:H0}

As the behaviour of $T_N$ does not depend on $\sigma^2$, we assume
that $\sigma^2 =1$; in particular, ${\bf \Sigma}= {\bf I}_K .$ Under
$H_0$, matrix $\hat {\bf R}$ is a complex Wishart matrix and it is
well-known (see for instance \cite{Meh04}) that the Jacobian of the
transformation between the entries of the matrix and the
eigenvalues/angles is given by the Vandermonde determinant
$\prod_{1\leq i<j\leq K}(x_j-x_i)^2.$ This yields the joint p.d.f. of
the ordered eigenvalues \eqref{eq:exactpdf0} where the normalizing
constant $Z(N,K,{\bf I}_K)$ is denoted by $Z_{K,N}^0$ for simplicity.

The celebrated result from Mar$\check{\textrm{c}}$enko and Pastur
\cite{MarPas67} states that the limit as $N,K\to \infty$ of
the c.d.f. $F_N(x)=\frac {\# \{i, \,\lambda_i\leq x\}}{K}$ associated
to the empirical distribution of the eigenvalues ($\lambda_i$) of
$\hat {\bf R}$ is equal to
$\Pmp\left((-\infty,x]\right)$ where
$\Pmp$ represents the
Mar$\check{\mathrm{c}}$enko-Pastur distribution:
\begin{equation}
 \label{eq:pmp}
\Pmp(dy) =
\mathbf{1}_{(\lambda^-,\lambda^+)}(y) \frac{\sqrt{(\lambda^+ -y)(y-\lambda^-)}}{2\pi c y}\,dy,
\end{equation}
with $\lambda^+=(1+\sqrt{c})^2$ and $\lambda^-=(1-\sqrt{c})^2$. This
convergence is very fast in the sense that the probability of
deviating from $\Pmp$ decreases as $e^{-N^2 × \textrm{const.}}.$ More
precisely, a simple application of the large deviations results in
\cite{BArGui97} yields that for any distance $d$ on the set of
probability measures on $\mathbb R$ compatible with the weak
convergence and for any $\delta >0,$
\begin{equation}
 \label{eq:gd0}
 \limsup_{N \to \infty} \frac{1}{N} \log \mathbb P_0\left(d(F_N, \Pmp) >\delta \right) = - \infty\ .
\end{equation}

Moreover, the largest eigenvalue $\lambda_1$ of $\hat {\bf R}$
converges a.s. to the right edge of the
Mar$\check{\mathrm{c}}$enko-Pastur distribution, that is $(1+\sqrt
c)^2.$ A further result due to Johnstone~\cite{Joh01} describes its
speed of convergence ($N^{-2/3}$) and its fluctuations (see also
\cite{Nad10pre} for complementary results). Let $\Lambda_1$ be defined
by:
\begin{equation}
  \Lambda_1 = N^{2/3}\left(\frac{\lambda_1 - (1+\sqrt{c_N})^2}{b_N}\right)\ ,
 \label{eq:lN}
\end{equation}
where $b_N$ is defined by
\begin{equation}
  \label{eq:bN}
  b_N := (1+\sqrt{c_N})\left(\frac{1}{\sqrt{c_N}}+1\right)^{1/3}\:,
\end{equation}
then $\Lambda_1$ converges in distribution toward a standard Tracy-Widom random variable with c.d.f.
$F_{TW}$ defined by:
\begin{equation}
  F_{TW}(x) =\exp\left( -\int_x^{\infty} (u-x)q^2(u)\,du\right)\, \quad \forall x\in \mathbb{R}\ ,
\label{eq:tw}
\end{equation}
where $q$ solves the Painlevé II differential equation:
$$
q''(x) = xq(x) + 2 q^3(x),\quad  q(x) \sim \textrm{Ai}(x) \quad \textrm{as} \quad x\to \infty
$$
and where Ai$(x)$ denotes the Airy function. In particular, $F_{TW}$
is continuous. The Tracy-Widom distribution was first introduced in
\cite{TraWid94,TraWid96} as the asymptotic distribution of the
centered and rescaled largest eigenvalue of a matrix from the Gaussian
Unitary Ensemble.

Tables of the Tracy-Widom law are available for instance in~\cite{bejan},
while a practical algorithm allowing to efficiently evaluate equation~(\ref{eq:tw})
can be found in~\cite{Bor09pre}.

\begin{remark} \label{rem:real1} In the case where the entries of
  matrix ${\bf Y}$ are real Gaussian random variables, the
  fluctuations of the largest eigenvalue are still described by a
  Tracy-Widom distribution whose definition slightly differs from the
  one given in the complex case (for details, see \cite{Joh01}).
\end{remark}

\subsection{Behaviour under hypothesis $H_1$}\label{subsec:H1}

In this case, the covariance matrix writes ${\bf \Sigma}=  \sigma^2 {\bf I}_K + {\bf h}{\bf
  h}^*$ and matrix ${\bf \hat R}$ follows a
\emph{single spiked} model. Since the behaviour of $T_N$ is not affected if the entries of ${\bf
  Y}$ are multiplied by a given constant, we find it convenient to
consider the model where ${\bf \Sigma}= {\bf I}_K + \frac{{\bf h}{\bf
  h}^*}{\sigma^2}$. Denote by
$$
\rho_K = \frac{\|{\bf h}\|^2}{\sigma^2}
$$
the {\em signal-to-noise} ratio (SNR), then matrix ${\bf \Sigma}$ admits the decomposition
${\bf \Sigma}= {\bf U} {\bf D} {\bf U}^*$ where ${\bf U}$ is a unitary matrix and
${\bf D}= {\rm diag} \left(\rho_K,1, \ldots,1\right).$
With the same change of variables from the entries of the matrix to the eigenvalues/angles
with Jacobian $\prod_{1\leq i<j\leq K}(x_j-x_i)^2,$ the p.d.f. of the ordered eigenvalues writes:
 \begin{equation}
  \label{eq:exactpdf1}
  p_{K}^{1,N}(x_{1:K}) = \frac{{\bs 1}_{(x_1\geq\dots\geq x_K\geq 0)}}{Z_{K,N}^1}\prod_{1\leq i<j\leq K}(x_j-x_i)^2 \:\prod_{j=1}^Kx_j^{N-K}
e^{-Nx_j}
I_K\left(\frac{N}{K}{\bf B}_K, {\bf X}_K\right)
\end{equation}
where the normalizing constant $Z(N,K,{\bf I}_K + {\bf hh^*})$ is
denoted by $Z_{K,N}^1$ for simplicity, ${\bf X}_K$ is the diagonal
matrix with eigenvalues $(x_1,\dots,x_K),$ ${\bf B}_K$ is the $K× K$
diagonal matrix with eigenvalues $(\frac{\rho_K}{1+
  \rho_K},0,\dots,0)$, and for any real diagonal matrices ${\bf C}_K,
{\bf D}_K,$ the spherical integral $I_K({\bf C}_K,{\bf D}_K)$ is
defined as
\begin{equation}
\label{eq:defIZ}
  I_K({\bf C}_K,{\bf D}_K) = \int e^{K \tr ({\bf C}_K {\bf Q} {\bf D}_K {\bf Q}^H)} dm_K({\bf Q}),
\end{equation}
with $m_K$ the Haar measure on the unitary group of size $K$ (see
\cite[Chapter 3]{Mui82} for details).

Whereas this rank-one perturbation does not affect the asymptotic
behaviour of $F_N$ (the convergence toward
$\Pmp$ and the deviations of the
empirical measure given by \eqref{eq:gd0} still hold under $\mathbb
P_1$), the limiting behaviour of the largest eigenvalue $\lambda_1$
can change if the signal-to-noise ratio $\rho_K$ is large enough.
\begin{assumption} \label{ass:rho}
The following constant $\rho\in{\bb R}$ exists:
\begin{equation}
  \label{eq:rho}
  \rho = \lim_{K\to\infty} \frac{\|{\bs h}\|^2}{\sigma^2}\ \left( =  \lim_{K\to\infty} \rho_K \right)\:.
\end{equation}
\end{assumption}
We refer to $\rho$ as the limiting SNR. We also introduce
$$
\lspiked = (1+\rho)\left( 1 + \frac c\rho\right).
$$
Under hypothesis $H_1$, the largest eigenvalue
has the following asymptotic behaviour as $N,K$ go to infinity:
\begin{align}
  & \lambda_1 \xrightarrow[H_1]{a.s.} \left\{
    \begin{array}[h]{ll}
    \lspiked  & \textrm{if } \rho > \sqrt{c}\:, \RR
    \lambda^+  & \textrm{otherwise,}
    \end{array}\right.
\label{eq:limTNd}
\end{align}
see for instance~\cite{BaiSil06} for a proof of this result. Note in
particular that $\lspiked$ is strictly larger than the right edge of
the support $\lambda^+$ whenever $\rho > \sqrt{c}$. Otherwise stated,
if the perturbation is large enough, the largest eigenvalue converges
outside the support of Mar$\check{\mathrm{c}}$enko-Pastur distribution.

\subsection{Limiting behaviour of $T_N$ under $H_0$ and $H_1$}

Gathering the results recalled in Sections \ref{subsec:H0} and
\ref{subsec:H1}, we obtain the following:
\begin{prop} Let Assumption \ref{ass:rho} hold true and assume that $\rho>\sqrt{c}$, then:
$$
T_N \xrightarrow[H_0]{a.s.} (1+\sqrt{c})^2\quad \textrm{and}
\quad T_N \xrightarrow[H_1]{a.s.} (1+\rho)\left(1+\frac c{\rho}\right)\quad \textrm{as}\ N,K\rightarrow\infty.
$$ 
\end{prop}

\section{Asymptotic threshold and $p$-values}\label{sec:asympt}

\subsection{Computation of the asymptotic threshold and $p$-value}

In Theorem~\ref{the:pfa} below, we take advantage of the convergence
results of the largest eigenvalue of $\hat{\bf R}$ under $H_0$ in the
asymptotic regime $N,K\to \infty$ to express the threshold and the
$p$-value of interest in terms of Tracy-Widom quantiles. Recall that
$\bar{F}_{TW}= 1- F_{TW}$, that $c_N=\frac KN$, and that $b_N$ is
given by \eqref{eq:bN}.



\begin{theo}
\label{the:pfa}
Consider a fixed level $\alpha\in (0,1)$ and let $\gamma_N$ be the
threshold for which the power of test~(\ref{eq:TN}) is maximum,
\emph{i.e.}  $p_N(\gamma_N)=\alpha$ where $p_N$ is defined by
\eqref{eq:ccdfexpr}. Then:
  \begin{enumerate}
  \item The following convergence holds true:
$$
\zeta_N\ \stackrel{\triangle}=\ \frac{N^{2/3}}{b_N} \left(\gamma_N - (1+\sqrt{c_N})^2 \right)\  
\xrightarrow[N,K\rightarrow\infty]{}\ \bar F_{TW}^{-1}(\alpha) \ .
$$
  \item The PFA of the following test
    \begin{equation}
      \label{eq:simpTest}
      T_N \gtrlessHU (1+\sqrt{c_N})^2 + \frac{b_N}{N^{2/3}}\: \bar F_{TW}^{-1}(\alpha)
    \end{equation}
    converges to $\alpha$.
  \item The $p$-value $p_N(T_N)$ associated with the GLRT can be approximated by:
    \begin{equation}
      \label{eq:approxpval}
      \tilde p_N(T_N) = \bar F_{TW}\left(\frac{N^{2/3}(T_N - (1+\sqrt{c_N})^2)}{b_N}\right)
    \end{equation}
    in the sense that $p_N(T_N) - \tilde p_N(T_N)\to 0$.
  \end{enumerate}
\end{theo}
\begin{remark}
  Theorem~\ref{the:pfa} provides a simple approach to compute both the
  threshold and the $p$-values of the GLRT as the dimension $K$ of the
  observed time series and the number $N$ of snapshots are large: The
  threshold $\gamma_N$ associated with the level $\alpha$ can be
  approximated by the righthand side of~(\ref{eq:simpTest}).
  Similarly, equation~(\ref{eq:approxpval}) provides a convenient
  approximation for the $p$-value associated with one experiment.
  These approaches do not require the tedious computation of the exact
  complementary c.d.f.~(\ref{eq:ccdfexpr}) and, instead, only rely on
  tables of the c.d.f. $F_{TW}$, which can be found for instance
  in~\cite{bejan} along with more details on the computational aspects
  (note that function $F_{TW}$ does not depend on any of the problem's
  characteristic, and in particular not on $c$). This is of importance
  in real-time applications, such as cognitive radio for instance,
  where the users connected to the network must quickly decide for the
  presence/absence of a source.
\end{remark}

\begin{proof}[Proof of Theorem \ref{the:pfa}]
  Before proving the three points of the theorem, we first describe
  the fluctuations of $T_N$ under $H_0$ with the help of the results
  in Section \ref{subsec:H0}. Assume without loss of generality that
  $\sigma^2=1$, recall that $T_N= \frac{\lambda_1}{K^{-1} \tr \hat{\bf
      R}}$ and denote by:
\begin{equation}
  \tilde T_N = \frac{N^{2/3}(T_N - (1+\sqrt{c_N})^2)}{b_N}
\label{eq:tildeTN}
\end{equation}
the rescaled and centered version of the statistics $T_N$. A direct
application of Slutsky's lemma (see for instance \cite{vaart}) together with the fluctuations of
$\lambda_1$ as reminded in Section \ref{subsec:H0} yields that
$\tilde T_N$
converges in distribution to a standard Tracy-Widom random variable
with c.d.f. $F_{TW}$ which is continuous over $\mathbb{R}$. Denote by
$F_N$ the c.d.f. of $\tilde T_N$ under $H_0$, then a classical result,
sometimes called Polya's theorem (see for instance \cite{BicMil92}),
asserts that the convergence of $F_N$ towards $F_{TW}$ is uniform over
$\mathbb{R}$:
\begin{equation}\label{eq:uniform}
\sup_{x\in \mathbb{R}} | F_N(x) -F_{TW}(x)| \xrightarrow[N,K\to\infty]{} 0\ .
\end{equation}
We are now in position to prove the theorem.

The mere definition of $\zeta_N$ implies that
$\alpha=p_N(\gamma_N)={\bar F}_N(\zeta_N)$. Due to \eqref{eq:uniform},
$\bar{F}_{TW}(\zeta_N)\to \alpha$. As $F_{TW}$ has a continuous
inverse, the first point of the theorem is proved.

The second point is a direct consequence of the convergence of $F_N$
toward the Tracy-Widom distribution: The PFA of
test~(\ref{eq:simpTest}) can be written as: $ {\bb P}_0\left(\tilde
  T_N > \bar F_{TW}^{-1}(\alpha)\right)$ which readily converges to
$\alpha$.

The third point is a direct consequence of \eqref{eq:uniform}:
$
p_N(T_N) - \tilde p_N(T_N) = \bar F_N(\tilde T_N) - \bar F_{TW}(\tilde T_N)\to 0\ .
$
This completes the proof of Theorem~\ref{the:pfa}.

\end{proof}


\section{Asymptotic analysis of the power of the test}
\label{sec:power}

In this section, we provide an
asymptotic analysis of the power of the GLRT as $N,K\to\infty$. As the
power of the test goes exponentially to zero, its error exponent is
computed with the help of the large deviations associated to the
largest eigenvalue of matrix $\hat{\bf R}$. The error exponent and
error exponent curve are computed in Theorem~\ref{the:errexp}, Section
\ref{sec:notion}; the large deviations of interest are stated in
Section \ref{sec:ld}.  Finally Theorem~\ref{the:errexp} is proved in
Section~\ref{sec:proof-EE}.

\subsection{Error exponents and error exponent curve}
\label{sec:notion}

The most natural approach to characterize the performance of a test is
to evaluate its power or equivalently its miss probability
\emph{i.e.}, the probability under $H_1$ that the receiver
decides hypothesis $H_0$. For a given level $\alpha\in(0,1)$, the miss probability writes:
\begin{equation}
  \label{eq:miss}
  \beta_{N,T}(\alpha) = \inf_{\gamma} \left\{ \mathbb{P}_1\left(T_N<\gamma\right),\ \gamma\ \textrm{such that}\ \mathbb{P}_0\left(T_N>\gamma\right) \leq \alpha\right\}\ .
\end{equation}
Based on Section~\ref{sec:exact}, the infimum is achieved when the
threshold coincides with $\gamma=p_N^{-1}(\alpha)$; otherwise stated,
$ \beta_{N,T}(\alpha) = \mathbb{P}_1\left(T_N <
  p_N^{-1}(\alpha)\right) $ (notice that the miss probability
  depends on the unknown parameters $\bs h$ and $\sigma^2$).  As
$\beta_{N,T}(\alpha)$ has no simple expression in the general case, we
again study its asymptotic behaviour in the asymptotic regime of
interest~(\ref{eq:regime}). It follows from Theorem~\ref{the:pfa} that
$p_N^{-1}(\alpha)\to \lambda^+=(1+\sqrt c)^2$ for $\alpha \in
(0,1)$. On the other hand, under hypothesis $H_1,$ $T_N$ converges
a.s. to $\lspiked$ which is strictly greater than $\lambda^+$ when the
ratio $\frac{\|{\bf h}\|^2}{\sigma^2}$ is large enough. In this case,
$\mathbb{P}_1\left(T_N < p_N^{-1}(\alpha)\right)$ goes to zero as it
expresses the probability that $T_N$ deviates from its limit
$\lspiked$; moreover, one can prove that the convergence to zero is
exponential in $N$:
\begin{equation}
\label{roughldp}
\mathbb{P}_1\left(T_N < x\right) \propto e^{-N I_{\rho}^+(x)}\qquad
\textrm{for}\quad  x\le \lspiked\ , 
\end{equation}
where $I_{\rho}^+$ is the so-called rate function associated to $T_N$.
This observation naturally yields the following definition of the
error exponent ${\mathcal E}_T$:
\begin{eqnarray}
\label{eq:defET}
{\cal E}_{T} &=& \lim_{N,K\to\infty} -\frac 1N \log \beta_{N,T}(\alpha) 
\end{eqnarray}
the existence of which is established in Theorem~\ref{the:errexp} below (as $N,K\to\infty$).
Also proved is the fact that ${\mathcal E}_T$ does not depend on $\alpha$.

The error exponent ${\cal E}_{T}$ gives crucial information on the
performance of the test $T_N$, provided that the level $\alpha$ is
kept fixed when $N,K$ go to infinity. Its existence strongly relies
on the study of the large deviations associated to the statistics
$T_N$.

In practice however, one may as well take benefit from the increasing
number of data not only to decrease the miss probability, but to
decrease the PFA as well. As a consequence, it is of practical
interest to analyze the detection performance when both the miss
probability and the PFA go to zero at exponential speed.  A couple
$(a,b)\in (0,\infty)× (0,\infty)$ is said to be an
\emph{achievable} pair of error exponents for the test $T_N$ if there
exists a sequence of levels $\alpha_N$ such that, in the asymptotic
regime~(\ref{eq:regime}),
\begin{equation}
  \lim_{N,K\to\infty} -\frac 1N \log \alpha_N = a \quad \textrm{and}\quad
  \lim_{N,K\to\infty} -\frac 1N \log \beta_{N,T}(\alpha_N) = b\ . \label{eq:exp}
\end{equation}
We denote by ${\cal S}_{T}$ the set of achievable pairs of error
exponents for test $T_N$ as $N,K\to\infty$. We refer to ${\cal
  S}_{T}$ as the \emph{error exponent curve} of $T_N$.

The following notations are needed in order to describe the error
exponent ${\mathcal E}_T$ and error exponent curve ${\cal S}_{T}$.

\begin{equation}
 \label{eq:Fp}
\left\{\begin{array}{lcll}
\mathbf{f}(x)&=& \int
\frac{1}{y-x}\Pmp(dy)& \textrm{for}\ x\in
\mathbb{R} \setminus (\lambda^-, \lambda^+)\\
\mathbf{F}^+(x) &=&\int \log(x-y)
\mathbb{P}_{\check{\mathrm{M}}\mathrm{P}}(dy)& \textrm{for}\ x \geq \lambda^+ 
\end{array}\right.\ .
\end{equation}

\begin{remark} Function $\mathbf{f}$ is the well-known Stieltjes
  transform associated to Mar$\check{\mathrm{c}}$enko-Pastur
  distribution and admits a closed-form representation formula. So
  does function $\mathbf{F}^+$, although this fact is perhaps less known.
  These results are gathered in Appendix \ref{appendix:representation}.
\end{remark}
Denote by $\Delta(\,\cdot \mid A)$ the convex indicator function \emph{i.e.} the
function equal to zero for $x\in A$ and to infinity
otherwise.
For  $\rho>\sqrt{c}$, define the function:
\begin{eqnarray}
\label{eq:Ip}
  I^+_\rho(x)   &= &\frac {x-\lspiked}{(1+\rho)}  -\left(1-c\right)\log\left(\frac x\lspiked \right)
- c\left(\mathbf{F}^+(x)-\F^+(\lspiked)\right)  +\Delta(x \mid [\lambda^+,\infty))\ .
\end{eqnarray}
Also define the function:
\begin{equation}
  \label{eq:IO}
  I^+_0(x)= x-\lambda^+  -\left(1-c\right)\log\left(\frac x{\lambda^+} \right)
- 2c\left(\mathbf{F}^+(x)-\F^+(\lambda^+)\right)  +\Delta(x \mid [\lambda^+,\infty))\ .
\end{equation}
We are now in position to state the main theorem of the section:
\begin{theo}
\label{the:errexp}
Let Assumption \ref{ass:rho} hold true, then:
\begin{enumerate}
\item For any fixed level $\alpha\in (0,1)$, the limit ${\cal E}_{T}$
  in~(\ref{eq:defET}) exists as $N,K\to \infty$ and satisfies:
  \begin{equation}
    \label{eq:errexp}
    {\mathcal E}_{T}= I^+_{\rho }(\lambda^+)
  \end{equation}
if $\rho>\sqrt{c}$ and ${\mathcal E}_{T}=0$ otherwise.

\item The error exponent curve of test $T_N$ is given by:
  \begin{eqnarray}
    \label{eq:SU}
    {\cal S}_{T} &=& \left\{ (I_0^+(x), I_\rho^+(x))\: : \: x \in (\lambda^+,\lspiked)\right\}
  \end{eqnarray}
  if $\rho>\sqrt{c}$ and ${\cal S}_{T}=\emptyset$ otherwise.
\end{enumerate}
\end{theo}

The proof of Theorem~\ref{the:errexp} heavily relies on the large
deviations of $T_N$ and is postponed to Section \ref{sec:proof-EE}.
Before providing the proof, it is worth making the following remarks.
\begin{remark}
  Several variants of the GLRT have been proposed in the literature,
  and typically consist in replacing the denominator $\frac
  1K\mathrm{tr}\,{\hat{\bf R}}$ (which converges toward $\sigma^2$) by a more
  involved estimate of $\sigma^2$ in order to decrease the
  bias~\cite{MesIT08,MesSP08,KRNA08,KRNA09}.  However, it can be
  established that the error exponents of the above variants are as
  well given by~(\ref{eq:errexp}) and~(\ref{eq:SU}) in the asymptotic
  regime.
\end{remark}
\begin{remark}
 The error exponent ${\mathcal E}_{T}$ yields a simple approximation of the miss probability
in the sense that $\beta_{N,T}(\alpha) \simeq e^{-N\,{\mathcal E}_{T}}$ as $N\to\infty$. It depends on
the limiting ratio $c$ and on the value of the SNR $\rho$ through the constant $\lspiked$.
In the high SNR case, the error exponent turns out to have a simple expression as a function of $\rho$.
If $\rho\to\infty$ then $\lspiked$ tends to infinity as well, which simplifies the expression
of rate function $I^+_{\rho }$. Using ${\bf F}^+(\lspiked) = \log\lspiked +o_\rho(1)$ where
$o_\rho(1)$ stands for a term which converges to zero as $\rho\to\infty$, it is straightforward to show that for
each $x\geq \lambda^+$,
$
I^+_{\rho }(x) = \log\rho -1-(1-c)\log x - c{\bf F}^+(x) + o_\rho(1)
$.
After some algebra, we finally obtain: 
$$
{\mathcal E}_{T} = \log\rho - (1+\sqrt c) - (1-c)\log(1+\sqrt c) - c\log\sqrt c + o_\rho(1)\ .
$$
At high SNR, this yields the following convenient approximation of the miss probability:
\begin{equation}
\beta_{N,T}(\alpha) \simeq \left( \psi(c)\, \rho \right)^N\ ,
\end{equation}
where $\psi(c) = e^{-(1+\sqrt c)}(1+\sqrt c)^{c-1}{c}^{-\frac c2}$.
\end{remark}

\subsection{Large Deviations associated to $T_N$}
\label{sec:ld}

In order to express the error exponents of interest, a rigorous
formalization of \eqref{roughldp} is needed.  Let us recall the
definition of a Large Deviation Principle: A sequence of random
variables $(X_N)_{N \in \mathbb N}$ satisfies a Large Deviation
Principle (LDP) under $\mathbb{P}$ in the scale $N$ with good rate
function $I$ if the following properties hold true:
\begin{itemize}
\item $I$ is a nonnegative function with compact level sets, i.e.
  $\{x, I(x) \leq t\}$ is compact for $t \in \mathbb R,$
\item for any closed set $F \subset \mathbb R,$ the following upper bound holds true:
\begin{equation}
 \label{eq:upperbdTN}
\limsup_{N \to \infty} \frac{1}{N} \log \mathbb{P}(X_N \in F)  \leq -\inf_F I\ .
\end{equation}
\item for any open set $G \subset \mathbb R,$ the following lower bound holds true:
\begin{equation}
 \label{eq:lowerbdTN}
\liminf_{N \to \infty} \frac{1}{N} \log \mathbb{P}(X_N \in G)  \geq -\inf_G I\ .
\end{equation}
\end{itemize}

For instance, if $A$ is a set such that $\inf_{\textrm{int}(A)} I =
\inf_{\textrm{cl}(A)} I (= \inf_A I)$, (where $\textrm{int}(A)$ and
$\textrm{cl}(A)$ respectively denote the interior and the closure of
$A$), then \eqref{eq:upperbdTN} and \eqref{eq:lowerbdTN} yield
\begin{equation}\label{eq:ldp-informal}
\lim_{N\to \infty} N^{-1} \log \mathbb{P}(X_N \in  A ) = - \inf_A I\ .
\end{equation}
Informally stated, 
$$
\mathbb{P} (X_N \in A ) \quad \propto\quad  e^{-N\inf_A I}\qquad \textrm{as}\ N\to \infty\ . 
$$
If, moreover $\inf_A I >0$ (which typically happens if the limit of
$X_N$ -if existing- does not belong to $A$), then probability
$\mathbb{P} (X_N \in A )$ goes to zero exponentially fast, hence a
{\em large deviation} (LD); and the event $\{X_N \in A\}$ can be referred
to as a {\em rare} event. We refer the reader to \cite{DemZei98} for
further details on the subject.



As already mentioned above, all the probabilities of interest 
are rare events as $N,K$ go to infinity related to large deviations
for $T_N.$ More precisely, Theorem~\ref{the:errexp} is merely a
consequence of the following Lemma.
\begin{lemma}
\label{lem:ldp}
Let Assumption \ref{ass:rho} hold true and let $N,K\to\infty$,
then:
\begin{enumerate}
\item Under $H_0,$ $T_N$ satisfies the LDP in the scale $N$ with
      good rate function $I_0^+$, which is increasing from 0 to
      $\infty$ on interval~$[\lambda^+,\infty)$.
    \item Under $H_1$ and if $\rho >\sqrt c$, $T_N$ satisfies the
      LDP in the scale $N$ with good rate function $I_\rho^+.$
      Function $I_\rho^+$ is decreasing from $I_\rho^+(\lambda^+)$ to
      0 on $[\lambda^+,\lspiked]$ and increasing from 0 to $\infty$ on
      $[\lspiked,\infty)$.
\item For any bounded sequence  $(\eta_N)_{N\geq 0}$,
\begin{equation}
 \lim_{N,K\to\infty} -\frac 1N \log {\bb P}_1\left( T_N < (1+\sqrt{c_N})^2 + \frac{\eta_N}{N^{2/3}} \right) = \left\{
    \begin{array}[h]{ll}
      I_\rho^+(\lambda^+) & \textrm{ if }\rho>\sqrt{c} \\
      0 & \textrm{otherwise.}
    \end{array}\right. \label{eq:ldpa}
\end{equation}
\item Let $x\in (\lambda^+,\infty)$ and let $(x_N)_{N\geq 0}$ be any
  real sequence which converges to $x$. If $\rho\leq\sqrt{c}$, then:
  \begin{equation} \label{eq:rhopetit}
   \lim_{N,K\to\infty} -\frac 1N \log {\bb P}_1\left( T_N < x_N \right) = 0\ .
      \end{equation}

\end{enumerate}
\end{lemma}
The proof of Lemma \ref{lem:ldp} is provided in
Appendix~\ref{app:ldp}.


\begin{remark}\label{rem:real2}
\begin{enumerate}
\item The proof of the large deviations for $T_N$ relies on the
  fact that the denominator $K^{-1}\tr\, \hat {\bf R}$ of $T_N$
  concentrates much faster than $\lambda_1$. Therefore, the large
  deviations of $T_N$ are driven by those of $\lambda_1$, a fact
  that is exploited in the proof.
\item In Appendix~\ref{app:ldp}, we rather focus on the large
  deviations of $\lambda_1$ under $H_1$ and skip the proof of Lemma
  \ref{lem:ldp}-(1), which is simpler and available (to some extent)
  in \cite[Theorem 2.6.6]{AGZ09}\footnote{see also the errata sheet for the
  sign error in the rate function on the authors webpage.}. Indeed, the proof of the LDP
  relies on the joint density of the eigenvalues. Under $H_1$, this
  joint density has an extra-term, the spherical integral, and is thus
  harder to analyze.
\item   Lemma \ref{lem:ldp}-(3) is not a mere consequence of Lemma
  \ref{lem:ldp}-(2) as it describes the deviations of $T_N$
  at the vicinity of a point of discontinuity of the rate function.
  The direct application of the LDP would provide a trivial
  lower bound ($-\infty$) in this case.
\item 
  In the case where the entries of matrix ${\bf Y}$ are real Gaussian
  random variables, the results stated in Lemma \ref{lem:ldp} will
  still hold true with minor modifications: The rate
  functions will be slightly different. Indeed, the computation of the
  rate functions relies on the joint density of the eigenvalues, which
  differs whether the entries of ${\bf Y}$ are real or complex.
\end{enumerate}
\end{remark}

\begin{figure}
  \centering
  \includegraphics[width=8cm]{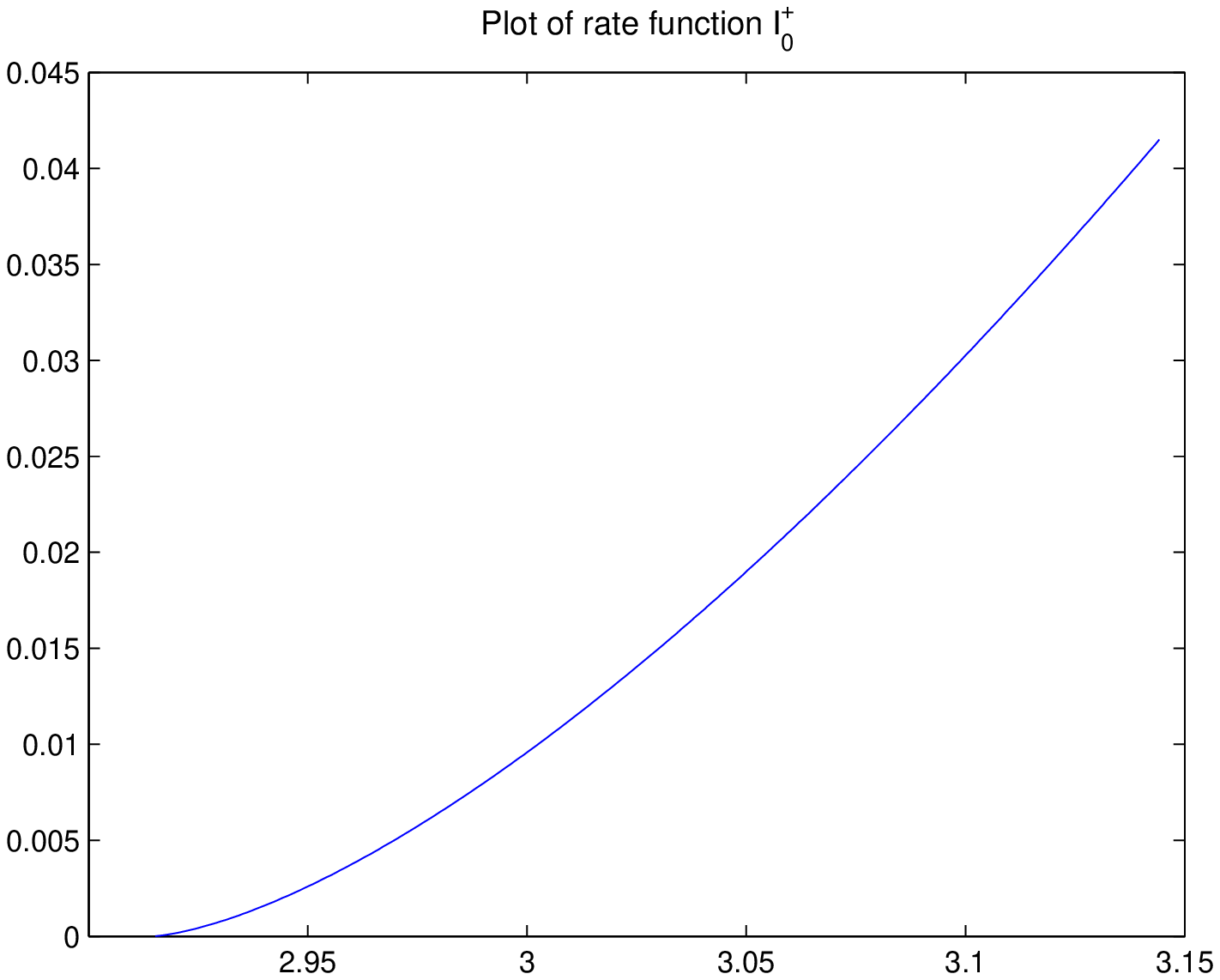}
  \includegraphics[width=8cm]{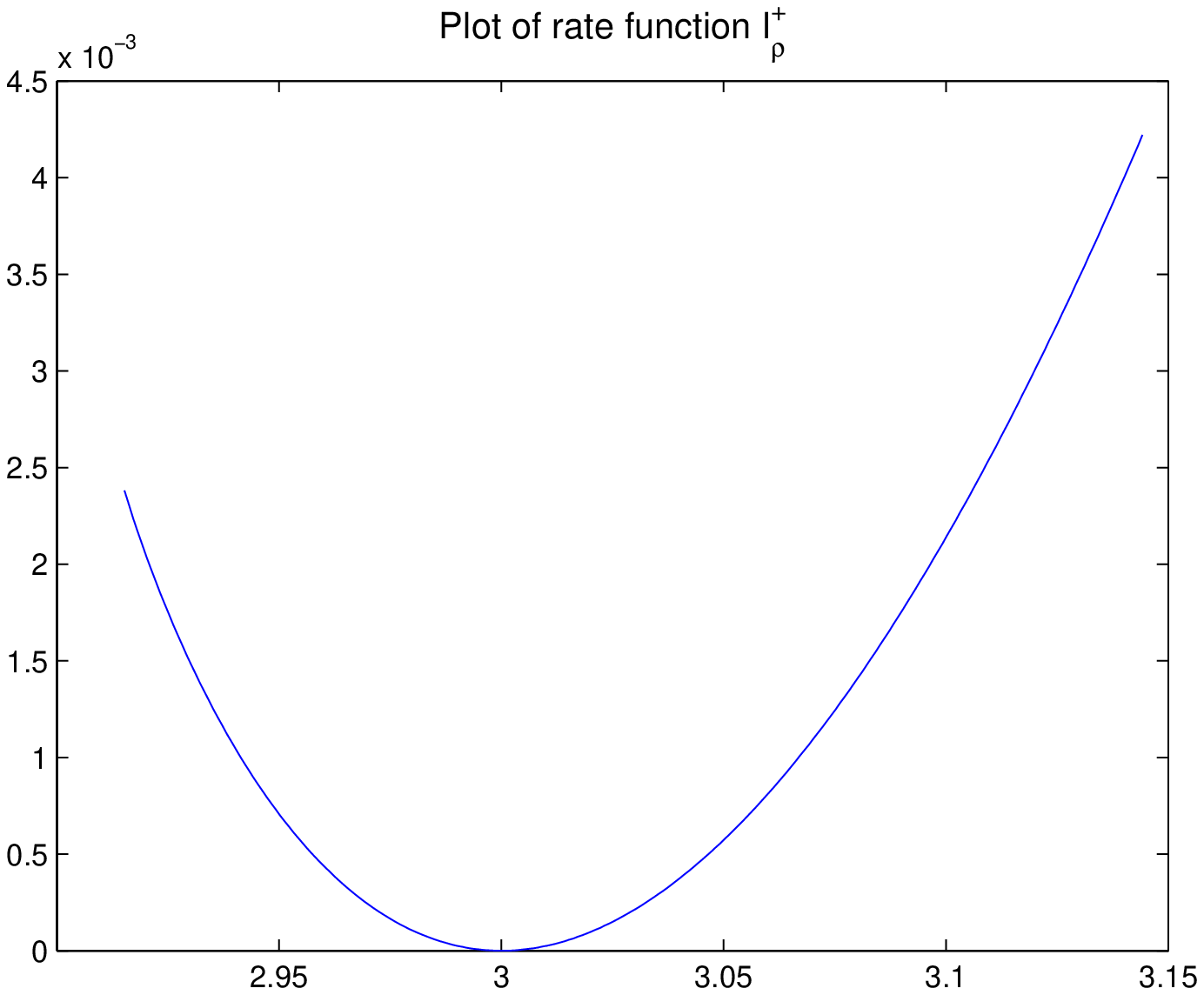}
  \caption{Plots of rate functions $I_0^+$ and $I_{\rho}^+$ in the
    case where $c=0.5$ and $\rho=1\,db$. In this case, $\lambda^+=2.9142$, $\lspiked=3$, $I_0^+(\lambda^+)=0$ and
$I_\rho^+(\lspiked)=0$.}
  \label{rate-functions}
 \end{figure}

\subsection{Proof of Theorem \ref{the:errexp}}\label{sec:proof-EE}
In order to prove (\ref{eq:errexp}), we must study the asymptotic
behaviour of the miss probability $\beta_{N,T}(\alpha) =
\mathbb{P}_1\left(T_N < p_N^{-1}(\alpha)\right)$ as $N,K\to
\infty$. Using Theorem~\ref{the:pfa}-(1), we recall that
\begin{equation}
\label{eq:betacloseLp}
\beta_{N,T}(\alpha) = \mathbb{P}_1\left(T_N < (1+\sqrt{c_N})^2 + \frac{\eta_N}{N^{2/3}}\right)
\end{equation}
where $c_N=\frac KN$ converges to $c$ and where $\eta_N$ is a deterministic sequence such that
$$
\lim_{N,K\to\infty} \eta_N = (1+\sqrt{c})\left(\frac{1}{\sqrt{c}}+1\right)^{1/3}\bar F_{TW}^{-1}(\alpha)\:.
$$
Hence, Lemma \ref{lem:ldp}-(3) yields the first point of
Theorem~\ref{the:errexp}.  We now prove the second point. Assume that
$\rho > \sqrt{c}$. Consider any $x\in(\lambda^+,\lspiked)$ and for
every $N,K$, consider the test function which rejects the null hypothesis
when $T_N>x,$
 \begin{equation}
   T_N \gtrlessHU x\ .
 \end{equation}
 Denote by $\alpha_N = {\bb P}_0(T_N>x)$ the PFA associated with this
 test. By Lemma \ref{lem:ldp}-(1) together with the continuity of the
 rate function at $x$, we obtain:
\begin{equation}
  \label{eq:proofa}
  \lim_{N,K\to\infty} -\frac 1N \log \alpha_N = \inf_{y \in [x, \infty)} I_0^+(y) = I_0^+(x)\ .
\end{equation}
The miss probability of this test is given by $\beta_{N,T}(\alpha_N) = {\bb P}_1(T_N<x)$. By Lemma \ref{lem:ldp}-(2),
\begin{equation}
  \label{eq:proofb}
  \lim_{N,K\to\infty} -\frac 1N \log \beta_{N,T}(\alpha_N) = \inf_{y \in (-\infty,x]} I_\rho^+(y) = I_\rho^+(x)\ .
\end{equation}
Equations~(\ref{eq:proofa}) and~(\ref{eq:proofb}) prove that
$(I_0^+(x), I_\rho^+(x))$ is an achievable pair of error exponents.
Therefore, the set in the righthand side of~(\ref{eq:SU}) is included
in ${\cal S}_{T}$. We now prove the converse. Assume that $(a,b)$ is
an achievable pair of error exponents and let $\alpha_N$ be a sequence such
that~(\ref{eq:exp}) holds. Denote by
$\gamma_N=p_N^{-1}(\alpha_N)$ the threshold associated with level $\alpha_N$. As
$I_0^+(x)$ is continuous and increasing from 0 to $\infty$ on interval
$(\lambda^+,\infty)$, there exists a (unique) $x\in(\lambda^+,\infty)$ such that $a=I_0^+(x)$.
We now prove that $\gamma_N$ converges to $x$ as $N$ tends to
infinity.  Consider a subsequence $\gamma_{\varphi(N)}$ which
converges to a limit $\gamma\in{\bb R}\cup\{\infty\}$. Assume that
$\gamma>x$. Then there exists $\epsilon>0$ such that
$\gamma_{\varphi(N)}>x+\epsilon$ for large $N$. This yields:
\begin{equation}
-\frac {1}{\varphi(N)} \log{\bb P}_0\left(T_{\varphi(N)} > \gamma_{\varphi(N)}\right)
\geq -\frac {1}{\varphi(N)} \log{\bb P}_0\left(T_{\varphi(N)} > x+\epsilon\right) \ .
\label{eq:extr}
\end{equation}
Taking the limit in both terms yields $I_0^+(x)\geq I_0^+(x+\epsilon)$ by
Lemma~\ref{lem:ldp}, which contradicts the fact that $I_0^+$ is
an increasing function. Now assume that $\gamma<x$.  Similarly,
\begin{equation}
-\frac {1}{\varphi(N)} \log{\bb P}_0\left(T_{\varphi(N)} > \gamma_{\varphi(N)}\right)
\leq -\frac {1}{\varphi(N)} \log{\bb P}_0\left(T_{\varphi(N)} > x-\epsilon\right)
\end{equation}
for a certain $\epsilon$ and for $N$ large enough. Taking the limit of both
terms, we obtain $I_0^+(x)\leq I_0^+(x-\epsilon)$ which leads to the same
contradiction. This proves that $\lim_N\gamma_N=x$.
Recall that by definition~(\ref{eq:exp}),
$$
  b=\lim_{N,K\to\infty} -\frac {1}{N} \log{\bb P}_1\left(T_{N} < \gamma_{N}\right)\ .
$$
As $\gamma_{N}$ tends to $x$, Lemma~\ref{lem:ldp} implies that the
righthand side of the above equation is equal to $I_\rho^+(x)>0$ if
$x\in(\lambda^+,\lspiked)$ and $\rho>\sqrt{c}$.  It is equal to 0 if $x\geq
\lspiked$ or $\rho\leq\sqrt{c}$.  Now $b>0$ by definition, therefore both
conditions $x\in(\lambda^+,\lspiked)$ and $\rho>\sqrt{c}$ hold.  As a conclusion,
if $(a,b)$ is an achievable pair of error exponents, then
$(a,b)=(I_0^+(x),I_\rho^+(x))$ for a certain $x\in(\lambda^+,\lspiked)$, and
furthermore $\rho>\sqrt{c}$.  This completes the proof of the second
point of Theorem~\ref{the:errexp}.

\section{Comparison with the test based on the condition number}
\label{sec:condition}

This section is devoted to the study of the asymptotic performances of
the test $U_N = \frac{\lambda_1}{\lambda_K}$, which is popular in
cognitive radio~\cite{zengliang,art-cardoso-2008,garello}. The main
result of the section is Theorem \ref{the:U}, where it is proved that
the test based on $T_N$ asymptotically outperforms the one based on
$U_N$ in terms of error exponent curves.

\subsection{Description of the test}

A different approach which has been introduced in several papers devoted to cognitive
radio contexts consists in rejecting the null hypothesis
for large values of the statistics $U_N$ defined by:
\begin{equation}
\label{eq:UN}
U_N = \frac{\lambda_1}{\lambda_K}\:,
\end{equation}
which is the ratio between the largest and the smallest eigenvalues of
$\hat {\bf R}.$ Random variable $U_N$ is the so-called \emph{condition
  number} of the sampled covariance matrix $\hat{\bf R}$.  As for
$T_N$, an important feature of the statistics $U_N$ is that its law
does not depend of the unknown parameter $\sigma$ which is the level
of the noise. Under hypothesis $H_0$, recall that the spectral measure
of $\hat {\bf R}$ weakly converges to the
Mar$\check{\mathrm{c}}$enko-Pastur distribution \eqref{eq:pmp} with
support $(\lambda^-,\lambda^+)$.  In addition to the fact that
$\lambda_1$ converges toward $\lambda^+$ under $H_0$ and $\lspiked$
under $H_1$, the following result related to the convergence of the
lowest eigenvalue is of importance (see for instance
\cite{yin-bai-kri-1988,bai-yin-1993}, \cite{BaiSil06}):
\begin{equation}
  \lambda_K \xrightarrow[]{a.s.}
    \lambda^-=\sigma^2(1-\sqrt{c})^2\
\label{eq:limUNd}
\end{equation}
under both hypotheses $H_0$ and $H_1$. Therefore, the statistics $U_N$ admits the following limits:
\begin{equation}
  U_N \xrightarrow[H_0]{a.s.} \frac{\lambda^+}{\lambda^-} = \frac{(1+\sqrt{c})^2}{(1-\sqrt{c})^2},
  \quad \textrm{and}\quad U_N \xrightarrow[H_1]{a.s.} \frac{\lspiked}{\lambda^-} \quad \textrm{for}\ \rho > \sqrt{c}\ .
\end{equation}
The test is based on the observation that the limit of $U_N$ under the
alternative $H_1$ is strictly larger than the ratio
$\lambda^+/\lambda^-$, at least when the SNR $\rho$ is large enough.

\subsection{A few remarks related to the determination of the
  threshold for the test $U_N$}

The determination of the threshold for the test $U_N$ relies on the
asymptotic independence of $\lambda_1$ and $\lambda_K$ under $H_0$.
As we shall prove below that test $U_N$ is asymptotically outperformed by test $T_N$, such a
study, rather involved, seems beyond the scope of this article.  For
the sake of completeness however, we describe unformally how to set
the threshold for $U_N$. Recall the definition of $\Lambda_1$ in
\eqref{eq:lN} and let $\Lambda_K$ be defined as:
$$
\Lambda_K =  N^{2/3}\frac{\left(\lambda_K - (1-\sqrt{c_N})^2\right)}{\left(\sqrt{c_N} -1\right)
\left( c_N^{-1/2} -1\right)^{1/3}}\ .
$$
Then both $\Lambda_1$ and $\Lambda_K$ converge toward Tracy-Widom
random variables. Moreover, 
$$
(\Lambda_1,\Lambda_K)\xrightarrow[N,K\rightarrow\infty]{} (X,Y)\ ,
$$
where $X$ and $Y$ are independent random variables, both distributed
according to $F_{TW}$\footnote{Such an asymptotic independence is
  not formally proved yet for $\bf{\hat R}$ under $H_0$, but is likely
  to be true as a similar result has been established in the case of
  the Gaussian Unitary Ensemble \cite{BiaDebNaj09},\cite{Bor09pre}.}.



As a corollary of the previous convergence, a direct application of
the Delta method \cite[Chapter 3]{Van98} yields the following convergence in distribution:
$$
N^{2/3} \left( 
\frac {\lambda_1}{\lambda_K} - \frac{(1+\sqrt{c_N})^2}{(1-\sqrt{c_N})^2}
\right) \rightarrow (a X + b Y)  \ ,
$$
where
$$
a = \frac{(1+\sqrt{c})}{(1-\sqrt{c})^2}\left( \frac 1{\sqrt{c}} +1\right)^{1/3}
\quad \textrm{and}\quad 
b = \frac{(1+\sqrt{c})^2}{(\sqrt{c}-1)^3}\left( \frac 1{\sqrt{c}} -1\right)^{1/3}\ ,
$$
which enables one to set the threshold of the test, based on the
quantiles of the random variable $a X + b Y$. In particular, following
the same arguments as in Theorem \ref{the:pfa}-1), one can prove that the
optimal threshold (for some fixed $\alpha \in (0,1)$), defined by
$
\mathbb{P}_0(U_N > \gamma_N ) = \alpha\ ,
$
satisfies 
$$
\xi_N\ \stackrel{\triangle}= \ N^{2/3} \left( \gamma_N -
  \frac{(1+\sqrt{c_N})^2}{(1-\sqrt{c_N})^2}\right)
\xrightarrow[N,K\rightarrow\infty]{} \bar{F}^{-1}_{aX + bY}(\alpha)\ .
$$
In particular, $\xi_N$ is bounded as $N,K\rightarrow\infty$.


\subsection{Performance analysis and comparison with the GLRT}

We now provide the performance analysis of the above test based on the
condition number $U_N$ in terms of error exponents.  In accordance
with the definitions of section~\ref{sec:notion}, we define the miss
probability associated with test $U_N$ as $\beta_{N,U}(\alpha) = \inf_{\gamma}
\mathbb{P}_1\left(U_N<\gamma\right)$ for any level $\alpha\in (0,1)$,
where the infimum is taken w.r.t. all thresholds $\gamma$ such that
$\mathbb{P}_0\left(U_N>\gamma\right)\leq \alpha$. We denote by ${\cal
  E}_{U}$ the limit of sequence $-\frac 1N \log \beta_{N,U}(\alpha)$
(if it exists) in the asymptotic regime~(\ref{eq:regime}).  We denote
by ${\cal S}_{U}$ the error exponent curve associated with test $U_N$
\emph{i.e.}, the set of couples $(a,b)$ of positive numbers for which
$-\frac 1N \log \beta_{N,U}(\alpha_N) \to b$ for a certain
sequence~$\alpha_N$ which satisfies $-\frac 1N \log \alpha_N \to a$.

Theorem~\ref{the:U} below provides the error exponents associated with
test $U_N$.  As for $T_N,$ the performance of the test is expressed in
terms of the rate function of the LDPs for $U_N$ under $\mathbb P_0$
or $\mathbb P_1$. These rate functions combine the rate functions for
the largest eigenvalue $\lambda_1$, i.e. $I_\rho^+$ and $I_0^+$
defined in Section \ref{sec:ld}, together with the rate function
associated to the smallest eigenvalue, $I^-$, defined below. As we
shall see, the positive rank-one perturbation does not affect
$\lambda_K$ whose rate function remains the same under $H_0$ and
$H_1$.

We first define:
\begin{equation}\label{eq:Fm}
  \mathbf{F}^-(x) = \int \log(y-x) d\Pmp (y)\quad \textrm{for}\ x\leq \lambda^-\ .
\end{equation}
As for $\mathbf{F}^+$, function $\mathbf{F}^-$ also admits a
closed-form expression based on $\mathbf{f}$, the Stieltjes transform
of Mar$\check{\textrm{c}}$enko-Pastur distribution (see Appendix
\ref{appendix:representation} for details).


Now, define for each $x\in{\bb R}$:
\begin{equation}\label{eq:Imoins}
I^-(x) = x - \lambda^- -\left( {1-c}\right)\log\left(\frac x{\lambda^-} \right)
-2c\left( \F^-(x) -\F^-(\lambda^-) \right) +\Delta(x |(0,\lambda^-])^ .
\end{equation}
If $\lambda_1$ and $\lambda_K$ were independent random variables, the
contraction principle (see e.g. \cite{DemZei98}) would imply that the
following functions
\begin{eqnarray*}
  \Gamma_{\rho}(t)= \inf_{(x,y)}\left\{ I_{\rho}^+(x) + I^-(y)\ :\quad \frac xy =t\right\} \quad \textrm{and}
\quad
  \Gamma_{0}(t)= \inf_{(x,y)}\left\{ I_{0}^+(x) + I^-(y)\ :\quad \frac xy =t\right\}
\end{eqnarray*}
defined for each $t\geq 0$, are the rate functions associated with the
LDP governing $\lambda_1/\lambda_K$ under hypotheses $H_1$ and $H_0$
respectively. Of course, $\lambda_1$ and $\lambda_K$ are not
independent, and the contraction principle does not apply.  However, a
careful study of the p.d.f. $p_{K,N}^{0}$ and $p_{K,N}^{1}$ shows that
$\lambda_1$ and $\lambda_K$ behave as if they were asymptotically
independent, from a large deviation perspective:

\begin{lemma}
\label{lem:ldpU}
Let Assumption \ref{ass:rho} hold true and let $N,K\to\infty$,
then:
\begin{enumerate}
\item Under $H_0,$ $U_N$ satisfies the LDP in the scale $N$ with
      good rate function $\Gamma_0$.
    \item Under $H_1$ and if $\rho >\sqrt c$, $U_N$ satisfies the
      LDP in the scale $N$ with good rate function $\Gamma_\rho.$
\item For any bounded sequence  $(\eta_N)_{N\geq 0}$,
\begin{equation}
 \lim_{N,K\to\infty} -\frac 1N \log {\bb P}_1\left( U_N < \frac{(1+\sqrt{c_N})^2}{(1-\sqrt{c_N})^2} + \frac{\eta_N}{N^{2/3}} \right) = \left\{
    \begin{array}[h]{ll}
      \Gamma_\rho(\lambda^+) & \textrm{ if }\rho>\sqrt{c} \\
      0 & \textrm{otherwise.}
    \end{array}\right. \label{eq:ldpa}
\end{equation}
Moreover, $\Gamma_{\rho}(\lambda^+)= I_{\rho}^+(\lambda^+)$.
\item Let $x\in (\lambda^+,\infty)$ and let $(x_N)_{N\geq 0}$ be any
  real sequence which converges to $x$. If $\rho\leq\sqrt{c}$, then:
  \begin{equation} \label{eq:rhopetit}
   \lim_{N,K\to\infty} -\frac 1N \log {\bb P}_1\left( T_N < x_N \right) = 0
      \end{equation}

\end{enumerate}
\end{lemma}

\begin{remark}\label{rem:lem2}
  In the context of Lemma \ref{lem:ldp}, both quantities $\lambda_1$
  and $\lambda_K$ deviate at the same speed, to the contrary of
  statistics $T_N$ where the denominator concentrated much faster than
  the largest eigenvalue $\lambda_1$. Nevertheless, proof of Lemma
  \ref{lem:ldpU} is a slight extension of the proof of
  Lemma~\ref{lem:ldp}, based on the study of the joint deviations
  $(\lambda_1,\lambda_K)$, the proof of which can be performed
  similarly to the proof of the deviations of $\lambda_1$. Once the
  large deviations established for the couple $(\lambda_1,\lambda_K)$,
  it is a matter of routine to get the large deviations for the ratio
  $\lambda_1/\lambda_K$. A proof is outlined in Appendix
  \ref{app:ldp-U}.
\end{remark}

We now provide the main result of the section.
\begin{theo}
  \label{the:U} Let Assumption \ref{ass:rho} hold true, then:
\begin{enumerate}
\item  For any fixed level $\alpha\in (0,1)$ and for each $\rho$, the error exponent
  ${\cal E}_{U}$ exists and coincides with~${\cal E}_{T}$.

\item The error exponent curve of test $U_N$ is given by:
  \begin{eqnarray}
    {\cal S}_{U} &=& \left\{ (\Gamma_0(t), \Gamma_\rho(t))\: : \: t \in \left(\frac{\lambda^+}{\lambda^-},\frac{\lspiked}{\lambda^-}\right)\right\}
  \end{eqnarray}
  if $\rho>\sqrt{c}$ and ${\cal S}_{U}=\emptyset$ otherwise.

\item The error exponent curve ${\cal S}_{T}$ of test $T_N$ uniformly dominates
${\cal S}_{U}$ in the sense that for each $(a,b)\in{\cal S}_{U}$ there exits $b'>b$ such that
$(a,b')\in{\cal S}_{T}$.
\end{enumerate}
\end{theo}
\begin{proof}
  The proof of items (1) and (2) is merely bookkeeping from the proof of
  Theorem \ref{the:errexp} with Lemma \ref{lem:ldpU} at hand.

  Let us prove item (3). The key observation lies in the following two facts:
\begin{eqnarray}\label{eq:key-observation1}
\forall x\in (\lambda^+,\lspiked),\quad \Gamma_{\rho}\left(\frac x{\lambda^-} \right) &=& I^+_{\rho}(x)\ ,\\
\forall x\in (\lambda^+,\lspiked),\quad \Gamma_0\left(\frac x{\lambda^-} \right) &<& I^+_0(x)\ .
\label{eq:key-observation2}
\end{eqnarray}
Recall that
\begin{eqnarray*}
  \Gamma_{\rho}\left( \frac x{\lambda^-} \right)&=& \inf_{(u,v)}\left\{ I_{\rho}^+(u) + I^-(v)\ :\quad \frac uv
=\frac x{\lambda^-}\right\}\\
&\stackrel{(a)}{\leq}& I_{\rho}^+(x) + I^-(\lambda^-)\quad =\,  I_{\rho}^+(x),
\end{eqnarray*}
where $(a)$ follows from the fact that $I^-(\lambda^-)=0$ and by
taking $u=x,v=\lambda^-$. Assume that inequality $(a)$ is strict. Due
to the fact that $I_{\rho}^+$ is decreasing, the only way to decrease
the value of $I^+_{\rho}(u) + I^-(v)$ under the considered constraint
$\frac uv =\frac x{\lambda^-}$ is to find a couple $(u,v)$ with $u>x$,
but this cannot happen because this would enforce $v>\lambda^-$ so
that the constraint $\frac uv = \frac x{\lambda^-}$ remains fulfilled,
and this would end up with $I^-(v)=\infty$. Necessarily, $(a)$ is an
equality and \eqref{eq:key-observation1} holds true.

Let us now give a sketch of proof for \eqref{eq:key-observation2}.
Notice first that $\frac{dI^+_0}{du}\mid_{u=x}>0$ (which easily
follows from the fact that $I^+_0$ is increasing and differentiable)
while $\frac{dI^-}{dv}\mid_{v\nearrow \lambda^-}=0$. This equality
follows from the direct computation:
\begin{eqnarray*}
\lim_{x\nearrow\lambda^-} \frac{I^-(x)}{x-\lambda^-}&=& 1 - \frac{1-c}{\lambda^-} -2c \left.\frac{d{\bf F}^-}{dx}\right|_{x\nearrow
  \lambda^-}\\
&=& 1 - \frac{1+\sqrt{c}}{1-\sqrt{c}} +2c{\bf f}(\lambda^-)\quad =\ 0\ ,
\end{eqnarray*}
where the last equality follows from the fact that $\frac{d{\bf
    F}^-}{dx}= -{\bf f}$ together with the closed-form expression for
${\bf f}$ as given in Appendix \ref{appendix:representation}. As previously, write:
\begin{eqnarray*}
  \Gamma_0\left( \frac x{\lambda^-} \right)&=& \inf_{(u,v)}\left\{ I_0^+(u) + I^-(v)\ :\quad \frac uv
=\frac x{\lambda^-}\right\}\\
&\stackrel{(a)}{\leq}& I_0^+(x) + I^-(\lambda^-)\quad =\,  I_0^+(x).
\end{eqnarray*}
Consider now a small perturbation $u=x-\delta$ and the related
perturbation $v=\lambda^--\delta'$ so that the constraint $\frac uv=
\frac x{\lambda^-}$ remains fulfilled. Due to the values of the
derivatives of $I^+_0$ and $I^-$ at respective points $x$ and
$\lambda^-$, the decrease of $I^+_0(x-\delta)$ will be larger than the
increase of $I^-(\lambda^--\delta')$, and this will result in the fact that
\begin{eqnarray*}
  \Gamma_0\left( \frac x{\lambda^-} \right)
\ \leq\ I_0^+(x-\delta) + I^-(\lambda^-+\delta')\ <\  I_0^+(x)\ ,
\end{eqnarray*}
which is the desired result, which in turn yields \eqref{eq:key-observation2}.

We can now prove Theorem \ref{the:U}-(3). Let $(a,b)\in {\mathcal
  S}_U$ and $(a,b')\in {\mathcal S}_T$, we shall prove that $b<b'$.
Due to the mere definitions of the curves ${\mathcal S}_U$ and
${\mathcal S}_T$, there exist $x\in (\lambda^+,\lspiked)$ and $t\in
(\lambda^+/\lambda^-,\lspiked/\lambda^-)$ such that
$a=I_0^+(x)=\Gamma_0(t)$. Eq. \eqref{eq:key-observation2} yields that
$\frac x{\lambda^-}< t$. As $I_{\rho}^+$ is decreasing, we have
$$
b'\ =\ I^+_{\rho}(x)\ >\ I^+_{\rho}(t \lambda^-)\ =\ \Gamma_{\rho}(t)\ =\ b\ ,
$$
and the proof is completed.
\end{proof}

\begin{remark}
  Theorem~\ref{the:U}-(1) indicates that when the number of data
  increases, the powers of tests $T_N$ and $U_N$ both converge to one
  at the same exponential speed ${\cal E}_{U}={\cal E}_{T}$, provided
  that the level $\alpha$ is kept fixed. However, when the level goes
  to zero exponentially fast as a function of the number of snapshots,
  then the test based on $T_N$ outperforms $U_N$ in terms of error
  exponents: The power of $T_N$ converges to one faster than the power
  of $U_N$. Simulation results for $N,K$ fixed sustain this claim
  (cf. Figure \ref{comparison}). This proves that in the context of
  interest ($N,K\to \infty$), the GLRT approach should be prefered to
  the test $U_N$.
\end{remark}


\section{Numerical Results}
\label{sec:numerical}
In the following section, we analyze the performance of the proposed
tests in various scenarios. 

\begin{figure}
  \centering
  \includegraphics[width=12cm]{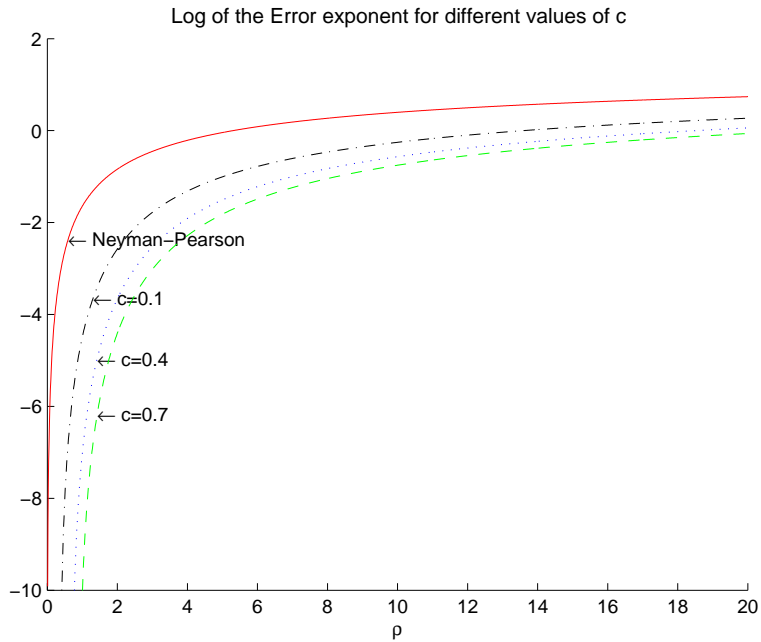}
  \caption{Computation of the logarithm of the error exponent
    ${\mathcal E}$ associated to the test $T_N$ for different values
    of $c$ (with ${\mathcal E}_{\rho}$ defined for $\rho\geq \sqrt{c}$
    and ${\mathcal E}_{\rho}\,|_{\rho=\sqrt{c}} =0$), and comparison
    with the optimal result (Neyman-Pearson) obtained in the case
    where all the parameters are perfectly known.}
  \label{errorexponent}
 \end{figure}

Figure \ref{errorexponent} compares the error exponent of test $T_N$
with the optimal NP test (assuming that all the parameters are known)
for various values of $c$ and $\rho$. The error exponent of the NP
test can be easily obtained using Stein's Lemma (see for
instance~\cite{Chen96}). 

\begin{figure}
  \centering
  \includegraphics[width=12cm]{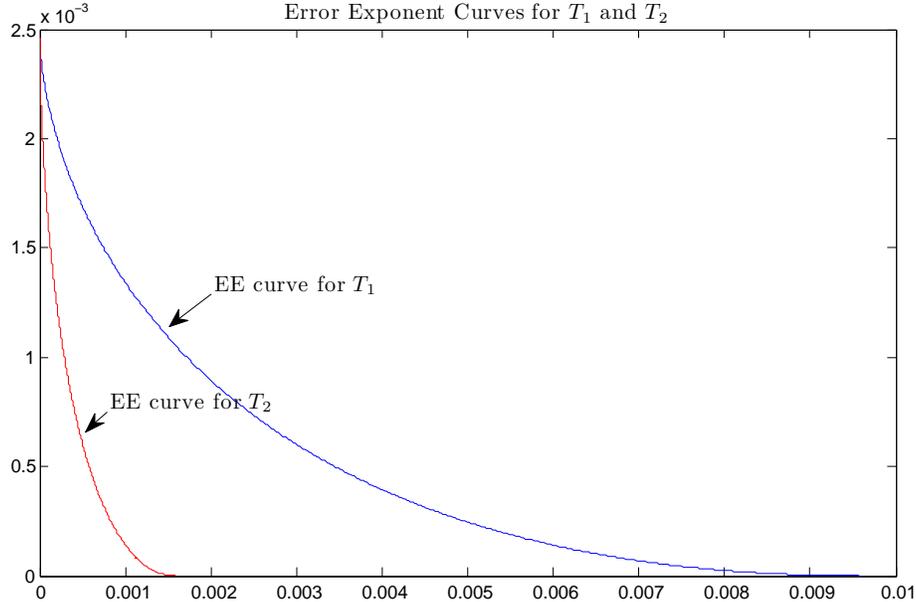} \caption{Error
    Exponent curves associated to the tests $T_N$ ($T_1$) and $U_N$
    ($T_2$) in the case where $c =\frac{1}{5}$ and $\rho = 10\, dB$.
    Each point of the curve corresponds to a given error exponent
    under $H_0$ ($X$ axis) and its counterpart error exponent under
    $H_1$ ($Y$ axis) as described in Theorem \ref{the:errexp}-(2) for $T_N$ and
    Theorem \ref{the:U}-(2) for $U_N$.}
  \label{ee-curve-T1-T2}
 \end{figure}

 In Figure \ref{ee-curve-T1-T2}, we compare the Error Exponent curves
 of both tests $T_N$ and $U_N$. The analytic expressions provided in
 \ref{the:errexp} and \ref{the:U} for the Error Exponent curves have
 been used to plot the curves. The asymptotic comparison clearly underlines the
 gain of using test $T_N$.

\begin{figure}
  \centering
  \includegraphics[width=12cm]{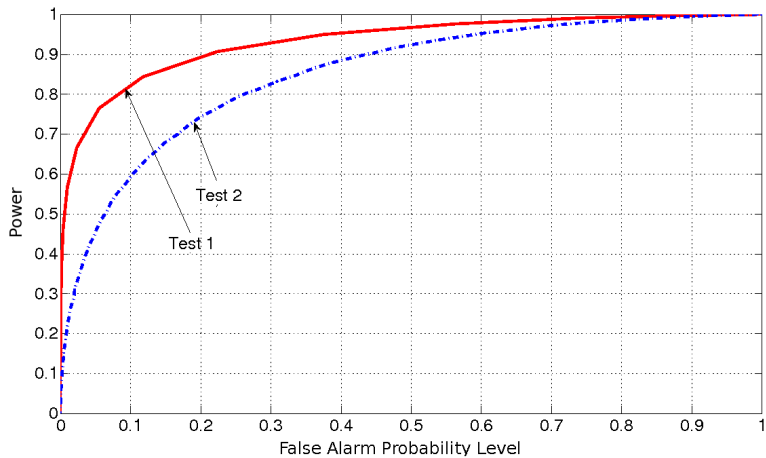}  
  \caption{Simulated ROC curves for $T_N$ (test 1) and $U_N$ (test 2)
    in the case where $K=10$, $N=50$ and $\rho=10\, dB$.}
  \label{comparison}
 \end{figure}

 Finally, we compare in Figure \ref{comparison} the powers (computed
 by Monte-Carlo methods) of tests $T_N$ and $U_N$ for finite values of
 $N$ and $K$. We consider the case where $K=10$, $N=50$ and $\rho=1$
 and plot the probability of error under $H_0$ versus the power of the
 test, that is $\alpha$ versus $\mathbb{P}_{1}(T_N \ge \gamma_N)$
 (resp. $\mathbb{P}_{1}(U_N \ge \gamma_N)$) where $\gamma_N$ is fixed by the following condition:
\begin{eqnarray*}
\mathbb{P}_{0}(T_N\geq \gamma_N)=\alpha \quad (\textrm{resp.} \  \mathbb{P}_{0}(U_N\geq \gamma_N)=\alpha)\ .
\end{eqnarray*}

\section{Conclusion}
\label{sec:conclusion}

In this contribution, we have analyzed in detail the GLRT in the case
where the noise variance and the channel are unknown. Unlike similar
contributions, we have focused our efforts on the analysis of the
error exponent by means of large random matrix theory and large
deviation techniques. Closed-form expressions were obtained and
enabled us to establish that the GLRT asymptotically outperforms the
test based on the condition number, a fact that is supported by
finite-dimension simulations. We also believe that the large
deviations techniques introduced here will be of interest for the
engineering community, beyond the problem addressed in this paper.

\section*{Acknowlegment}
We thank Olivier Cappé for many fruitful discussions related to the GLRT.

\appendices

\section{Proof of Lemma~\ref{lem:ldp}: Large deviations for $T_N$}
\label{app:ldp}

The large deviations of the largest eigenvalue of large random
matrices have already been investigated in various contexts, Gaussian
Orthogonal Ensemble \cite{BDG01} and deformed Gaussian ensembles
\cite{Mai07}.  As mentionned in \cite[Remark 1.2]{Mai07}, the proofs
of the latter can be extended to complex Wishart matrix models, that
is random matrices $\bf{\hat R}$ under $H_0$ or $H_1$. 

In both cases, the large deviations of $\lambda_1$ rely on a close
study of the density of the eigenvalues, either given by
\eqref{eq:exactpdf0} (under $H_0$) or by \eqref{eq:exactpdf1} for the
spiked model (under $H_1$). The study of the spiked model, as it
involves the study of the asymptotics of the spherical integral (see
Lemma \ref{IZlemma} below), is more difficult. We therefore focus on
the proof of the LDP under $H_1$ (Lemma~\ref{lem:ldp}-(2)) and omit
the proof of Lemma~\ref{lem:ldp}-(1). Once Lemma~\ref{lem:ldp}-(2) is
proved, proving Lemma~\ref{lem:ldp}-(1) is a matter of bookkeeping,
with the spherical integral removed at each step.

Recall that $\lambda_1\geq \cdots\geq \lambda_K$ are the ordered
eigenvalues of $\hat{\bf R}$ and that $T_N$ is the statistics defined
in \eqref{eq:mu}.

In the sequel, we shall prove the upper bound of the LDP in
Lemma~\ref{lem:ldp}-(2) (which gives also the upper bound in Lemma~\ref{lem:ldp}-(3)).
 The proof of the
lower bound in Lemma~\ref{lem:ldp}-(3) requires more precise arguments
than the lower bound of the LDP. One has indeed to study what
  happens at the vicinity of $\lambda^+,$ which is a point of
  discontinuity of the rate function $I_\rho^+.$ Thus, we skip the
proof of the lower bound of the LDP in Lemma~\ref{lem:ldp}-(2) to
avoid repetition. Note that the proof of Lemma~\ref{lem:ldp}-(4) is a
mere consequence of the fact that $T_N$ converges a.s. to $\lambda^+$
if $\rho \leq \sqrt c,$ thus $\mathbb P_1(T_N <x_N)$ converges to 1
whenever $x_N$ converges to $x > \lambda^+.$

For sake of simplicity and with no loss of generality as the law of
$T_N$ does not depend on $\sigma,$ we assume all along this appendix
that $\sigma^2 =1.$ We first recall important asymptotic results for
spherical integrals.

\subsection{Useful facts about spherical integrals}
Recall that the joint distributions of the ordered eigenvalues under
hypothesis $H_0$ and $H_1$ are respectively given by
\eqref{eq:exactpdf0} and \eqref{eq:exactpdf1}. In the latter, the
so-called spherical integral \eqref{eq:defIZ} is introduced.  We
recall here results from \cite{Mai07} related to the asymptotic
behaviour of the spherical integral in the case where
one diagonal matrix is of rank one and the other has the
limiting distribution $\Pmp$.
We first introduce the function defined for $x \geq \lambda^+$ by:
\begin{equation}\label{eq:jrho}
J_\rho(x) = \left\{
\begin{array}{ll}
 \frac{\rho}{c} - \log\left( \frac{\rho}{c(1+\rho)}\right) - {\bf F}^+(\lspiked),& \textrm{ if } \rho \leq \sqrt c \textrm{ and }
\lambda^+  \leq x \leq \lspiked, \\
\frac{\rho x}{c(1+\rho)} -1- \log\left( \frac{\rho}{c(1+\rho)}\right) - {\bf F}^+(x),& \textrm{ otherwise.}
\end{array}
\right.
\end{equation}

Consider a $K$-tuple $(x_1,\cdots, x_K)$ and denote by $\hat
\pi_{K,{\bf x}} = \frac{1}{K-1} \sum_{i=2}^N \delta_{x_2}$ the
empirical distribution associated to $(x_2,\cdots, x_K)$; let $d$ be
a metric compatible with the topology of weak convergence of measures
(for example the Dudley distance - see for instance \cite{Dud02}). A
strong version of the convergence of the spherical integral in the
exponential scale with speed $N$, established in \cite{Mai07} can be
summarized in the following Lemma:

\begin{lemma} \label{IZlemma}
Assume that $N,K\to \infty$ and $\frac KN \to c\in (0,1)$ and let Assumption \ref{ass:rho} hold true.
Let $x_1\geq x_2\geq \cdots \geq x_K\geq 0$ and $\delta >0$. If, for $N$ large enough, $| x_1-x|\leq \delta$
and $d(\hat \pi_{K,\bf x}, \Pmp) \leq N^{-1/4}$ then:
$$
\left| \frac{1}{N} \log I_K\left(\frac{N}{K}{\bf B}_K, {\bf X}_K\right) - c J_\rho(x)\right|\leq \delta,
$$
where $J_{\rho}$ is given by \eqref{eq:jrho}, ${\bf B}_K = {\rm
  diag}\left(\frac{\rho_K}{1+\rho_K}, 0, \ldots,0\right)$ and
${\bf X}_K=\mathrm{diag}(x_1,\cdots, x_K)$.
\end{lemma}

Recall that the spherical integral $I_K$, defined in \eqref{eq:defIZ},
appears in the joint density \eqref{eq:exactpdf1} of the eigenvalues
under $H_1$. Lemma \ref{IZlemma} provides a simple asymptotic
equivalent $cJ_{\rho}(x)$ of the normalized integral $N^{-1} \log
I_K$. Roughly speaking, this will enable us to replace $I_K$ by the
quantity $e^{-N\times cJ_{\rho}(x)}$ when establishing the large deviations of
$\lambda_1$, which rely on a careful study of density
\eqref{eq:exactpdf1}.


\subsection{Proof of Lemma \ref{lem:ldp}-(2)}

In order to establish the LDP under hypothesis $H_1$ and condition
$\rho > \sqrt c$,  (that is the
bounds \eqref{eq:upperbdTN} and \eqref{eq:lowerbdTN}), we first notice
that intervals $(x, x+\delta)$ for $x, \delta \in \mathbb R^+$ form a
basis of the topology of $\mathbb R^+$. The LDP will be therefore a
consequence of the following bounds:
\begin{itemize}
\item(Exponential tightness) there exists a function $f : \mathbb R^+
  \to \mathbb R^+$ going to infinity at infinity such that for
  all $N$,
\begin{equation} \label{exptight}
 \mathbb P_1\left( \lambda_1 \geq M\right) \ \leq \ e^{-N f(M)} \ .
\end{equation}
Condition \eqref{exptight} is technical (see for instance \cite[Lemma
1.2.18]{DemZei98}): Instead of proving the large deviation upper bound
for every closed set, the exponential tightness \eqref{exptight}, if
established, enables one to restrict to the compact sets.

\item(Upper bound) For any $x$, for any $M$ such that $0<x < M, $
\begin{equation} \label{upperbd}
\lim_{\delta \downarrow 0} \limsup_{N,K \to \infty}
\frac 1 N \log \mathbb P_1 \left(x \leq T_N \leq x+\delta, \,\lambda_1 \leq M\right)
\ \leq \ - I_\rho^+(x) \ ,
\end{equation}
Due to the exponential tightness, it is sufficient to establish the
upper bound for compact sets. As each compact can be covered by a
finite number of balls, it is therefore sufficient to establish upper
estimate \eqref{upperbd} in order to establish the LD upper bound.

\item(Lower bound)    For any $x$,
\begin{equation} \label{lowerbd}
\lim_{\delta \downarrow 0} \liminf_{N,K \to \infty}
\frac 1 N \log \mathbb P_1\left(x \leq T_N \leq x+\delta\right)
\ \geq\ - I_\rho^+(x)\ .
\end{equation}
The fact that \eqref{lowerbd} implies the LD lower bound
\eqref{eq:lowerbdTN} is standard in LD and can be found in
\cite[Chapter 1]{DemZei98} for instance.
\end{itemize}
As the arguments are very similar to the ones developed in
\cite{Mai07}, we only prove in detail the upper bound
\eqref{upperbd}. Proofs of \eqref{exptight} and \eqref{lowerbd} are
left to the reader.

The idea is that the
empirical measure $\hat \pi_{K,{\boldsymbol \lambda}}:= \frac{1}{K-1}
\sum_{j=2}^K \delta_{\lambda_j}$ (of all but the largest eigenvalues)
and the trace concentrate faster than the
largest eigenvalue. In the exponential scale with speed $N$, $\hat
\pi_{K,{\boldsymbol \lambda}}$ and the trace can be considered as equal to their
limit, respectively $P_{\check{\mathrm M}{\mathrm P}}$ and 1. In
particular, the deviations of $T_N$ arise from those of the largest
eigenvalue and they both satisfy the same LDP with the same rate
function $I_\rho^+$. We therefore isolate the terms depending on $\lambda_1$
and gather the others through their empirical measure $\hat \pi_{K,\boldsymbol \lambda} .$

Recall the notations introduced in \eqref{eq:exactpdf0} and
\eqref{eq:exactpdf1} and let $x > \lambda^+$, $\delta >0$. Consider the following domain:
$$
{\mathcal D}=\left\{ (x_1,\cdots, x_K)\in [0,M]^{K},\ \frac{K x_1}{x_1+\cdots +x_K} \in (x, x+\delta) \right\}
$$
For $N$
large enough:
\begin{eqnarray*}
 \lefteqn{\mathbb P_1(x \leq T_N \leq x+\delta,\,\, \lambda_1 \leq M )=
\int_{\mathcal D} dp_{K,N}^{1}(x_{1:K})}\\
&=& \frac{1}{Z_{K,N}^1}
\int_{\mathcal D} dx_1\,   
e^{-Nx_1}e^{(N-K) \log x_1}e^{2(K-1)\int \log(x_1-u)d\hat\pi_{K,{\bf x}}(u)} \\
&&× I_K\left(\frac{N}{K} {\bf B}_K, {\bf X}_K\right)
\prod_{1<i<j} |x_i-x_j|^2 e^{-N \sum_{j=2}^K x_j} \prod_{j=2}^K x_j^{N-K} d\,x_{2:K}× {\bs 1}_{(x_1\geq\dots\geq x_K\geq 0)} \\
&=&  \frac{\left(1- \frac{1}{N}\right)^{(K-1)(N-1)} Z_{K-1,N-1}^0 }{Z_{K,N}^1} \int_{\mathcal D}
dx_1
e^{-Nx_1}e^{(N-K) \log x_1}e^{2(K-1)\int \log(x_1-u)d\hat\pi_{K,{\bf y}}(u)} \\
&&× I_K\left(\frac{N}{K}{\bf B}_K,{\bf X}_K\right) dp_{K-1,N-1}^{0}(y_{2:K}),
\end{eqnarray*}
where we performed the change of variables $y_i := \frac{N}{N-1}x_i$
for $i=2:K$, and the related modifications $\hat \pi_{K,{\bf x}}\leftrightarrow \hat \pi_{K,{\bf y}}$
and ${\bf X}_K=\mathrm{diag}\left(x_1,\frac{N-1}N y_2,\cdots, \frac{N-1}N y_2\right)$.
Note also that strictly speaking, the domain of integration ${\mathcal D}$ would
express differently with the $y_i$'s and in particular, we should have
changed constant $M$ which majorizes the $x_i$'s into a larger
constant as the $y_i$'s can theoretically be slightly above $M$ - we
keep the same notation for the sake of simplicity.

To proceed, one has to study the asymptotic behaviour of the normalizing constant:
$$
\displaystyle \frac{\left(1- \frac{1}{N}\right)^{(K-1)(N-1)}
  Z_{K-1,N-1}^0 }{Z_{K,N}^1} \ ,
$$
which turns out to be difficult. Instead of establishing directly
the bounds \eqref{exptight}-\eqref{lowerbd}, we proceed as in
\cite{Mai07} and establish similar bounds replacing the probability
measures $\mathbb P_1$ by the measures $\mathbb Q_1$ defined as:
$$
\mathbb Q_1
:= \frac{Z_{K,N}^1}{Z_{K-1, N-1}^0 \left(1- \frac{1}{N}\right)^{(K-1)(N-1)}} \mathbb P_1
$$
and the rate function $I_\rho^+$ by the function $G_\rho$ defined by:
$$ G_\rho(x) = \frac{x}{1+\rho} -(1-c) \log x  - c {\bf F}^+(x) +c + c\log\left( \frac{\rho}{c(1+\rho)} \right)
$$
for $x > \lambda^+$. Notice that these positive measures $\mathbb Q_1$
are not probability measures any more, and as a consequence, the
function $G_\rho$ is not necessarily positive and its infimum might not be
equal to zero, as it is the case for a rate function.

Writing the upper bound for $\mathbb Q_1$, we obtain:
\begin{eqnarray*}
\lefteqn{ \mathbb Q_1(x \leq T_N \leq x+\delta,\,\, \lambda_1 \leq M )}\\
& \leq&
\int_{\mathcal D} dx_1 e^{-N \Phi(x_1,c_N, \hat \pi_{K,{\bf y}})}
I_K\left(\frac{N}{K}{\bf B}_K, {\bf X}_K\right) dp_{K-1,N-1}^{0}(y_{2:K}),
\end{eqnarray*}
where, for any compactly supported probability measure $\mu$ and any real number $y$ greater than the right edge
of the support of $\mu,$
$$\Phi(y,c_N,\mu) = -y +(1-c_N) \log y + 2c_N \int \log (y-\lambda) d\mu(\lambda).$$

Let us now localise the empirical measure $\hat \pi_{K,{\bf y}}$
around $\Pmp$\footnote{Notice that if $\hat \pi_{K,{\bf x}}$ is close to $\Pmp$, so is $\hat \pi_{K,{\bf y}}$
due to the change of variable $y_i=\frac N{N-1} x_i$.} and the trace around 1.
The continuity and convergence properties
 of the spherical integral
recalled in Lemma \ref{IZlemma} yield, for $K$ large enough:
\begin{eqnarray}
 \mathbb Q_1(x \leq T_N \leq x+\delta\ ,\ \lambda_1 \leq M )
&\leq& \int_x^{x+\delta} dx_1 \int_{\mathcal E}
e^{-N \Phi(x_1,c_N, \hat \pi_{K,{\bf y}})} e^{Nc(J_\rho(x_1)+ \delta)} dp_{K-1,N-1}^{0}(y_{2:K}) \nonumber\\
&&+ 4^K M^{N+K} e^{NM\frac{\rho_K}{1+\rho_K}} \int_{\mathcal E^C}  dp_{K-1,N-1}^{0}(y_{2:K}), \label{eq:upM}
\end{eqnarray}
with $$ \mathcal E := \left\{(y_2,\cdots, y_K) \in [0,M]^{K-1},\ d(\hat \pi_{K,{\bf y}},\Pmp)
\leq  \frac 1{N^{1/4}}
\quad\mathrm{and}\quad \frac{1}{K}\sum_{j=2}^K y_j \in \left[1-\delta^2, 1+\delta^2 \right]\right\}.$$
The second term in \eqref{eq:upM} is easily obtained considering the fact that all the eigenvalues are less than $M$
so that for $1 \leq j \leq K,$
$ |x_1-x_j| \leq 2M,$ $x_j^{N-K} \leq M^{N-K}$ and $(UX_KU^*)_{11}\leq M.$
Now, standard concentration results under $H_0$ yield that:
$$ \limsup_{N,K \to \infty} \frac{1}{N} \log   \mathbb P_0\left( \hat \pi_{K,{\bf \lambda}} \notin B(\Pmp, N^{-1/4})
\ \mathrm{ or }\  \frac{1}{K}\sum_{j=2}^K \lambda_j \notin \left[1-\delta^2, 1+\delta^2\right]\right) = -\infty.
$$
More precisely, one knows using \cite{GuiZei00} that the empirical measure   $\frac{1}{K}\sum_{j=2}^K \lambda_j$
is close enough to its expectation and then using \cite{Bai93b} one knows that the expectation is close enough
to its limit $\Pmp.$ The arguments are detailed in the Wigner case in \cite{Mai07} and we do not give more details here.\\

As $c_N\to c$ for $N,K \to \infty$, $c \mapsto
\Phi(y,c,\mu) $ is continuous and $\mu \mapsto \Phi(y,c,\mu)$ is lower
semi-continuous, we obtain:
$$
\limsup_{N,K \to \infty} \frac{1}{N} \log \mathbb Q_1(x \leq
\lambda_1 \leq x+\delta\ ,\ \lambda_1 \leq M )\ \leq\
\sup_{u\in [x , x+\delta]} \left(\Phi(u,c, P_{\check{\mathrm M}{\mathrm
      P}}) + c J_\rho\left( u\right)\right) + 2\delta.
$$
By continuity in $u$ of the two involved functions, we finally get:
$$
\lim_{\delta\downarrow 0}\limsup_{N,K \to \infty} \frac{1}{N}
\log \mathbb Q_1(x \leq \lambda_1 \leq x+\delta\ ,\ \lambda_1
 \leq M )\ \leq\ \Phi(x,c, \Pmp) + c
J_\rho\left( x\right) = G_\rho(x)\ ,
$$
and the counterpart of Eq. \eqref{upperbd} is proved for $\mathbb Q_1$
and function $G_\rho$.  The proof of the lower bound is quite similar
and left to the reader. It remains now to recover \eqref{upperbd}. As
$\mathbb P_1$ is a probability measure and the whole space $\mathbb
R^+$ is both open and closed, an application of the upper and lower
bounds for $\mathbb Q_1$ immediately yields:
\begin{eqnarray}
  \lefteqn{\liminf_{N,K  \to \infty} \frac{1}{N} 
\log\frac{Z_{K,N}^{1}}{Z_{K-1,N-1}^{0} \left(1- \frac{1}{N}\right)^{(K-1)(N-1)}}
    \mathbb P_1(T_N \in \mathbb R^+)}\nonumber \\
  &=& \limsup_{N,K  \to \infty} \frac{1}{N}  \log \frac{Z_{K,N}^{1}}{Z_{K-1,N-1}^{0}\left(1- \frac{1}{N}\right)^{(K-1)(N-1)}}
  \mathbb P_1(T_N \in \mathbb R^+)\nonumber \\
  &=& \lim_{N,K  \to \infty} \frac{1}{N} \log \frac{Z_{K,N}^{1}}{Z_{K-1,N-1}^{0} \left(1- \frac{1}{N}\right)^{(K-1)(N-1)}}
\nonumber\\
  &=& -\inf_{\mathbb R^+} G_\rho\ .\label{eq:conv-constante}
\end{eqnarray}
This implies that the LDP holds for $\mathbb P_1$ with rate function
$G_\rho -\inf_{\mathbb R^+} G_\rho$.

It remains to check that $I_\rho^+=G_\rho -\inf_{\mathbb R^+} G_\rho$,
which easily follows from the fact to be proved that:
\begin{equation}\label{eq:Ginf}
\inf_{x \in [\lambda^+, \infty)}G_\rho(x)
= G_\rho(\lspiked)\ .
\end{equation}
We therefore study the variations of $G_{\rho}$ over $[\lambda^+,
\infty)$.  Note that $({\bf F}^+)^\prime = -{\bf f}$, and thus that
$G_{\rho}'(x)= (1+\rho)^{-1} - (1-c)x^{-1} + c{\bf f}(x)$. Function
${\bf f}$ being a Stieltjes transform is increasing for $x >
\lambda^+$, and so is $G_{\rho}'$, whose limit at infinity is
$(1+\rho)^{-1}$. Straightforward but involved computations using the
explicit representation \eqref{eq:fplus} for ${\bf f}$
yield that $G'_{\rho}(\lspiked)=0$. Therefore, $G_{\rho}$ is
decreasing on $[\lambda^+,\lspiked]$ and increasing on
$[\lspiked,\infty)$, and \eqref{eq:Ginf} is proved.

This concludes the proof of the upper bound in Lemma
\ref{lem:ldp}-(2).  The proof of Lemma~\ref{lem:ldp}-(1) is very
similar and left to the reader.

\subsection{Proof of Lemma~\ref{lem:ldp}-(3)}

The proof of this point requires an extra argument as we study the
large deviations of $T_N$ near the point $(1+\sqrt{c})^2$ where the
rate function is not continuous. In particular, the limit
\eqref{eq:ldpa} does not follow from the LDP already established. As
we shall see when considering $\mathbb{P}_1\left(T_N <
  (1+\sqrt{c_N})^2 + \eta_N N^{-2/3} \right)$, the fact that the
scale $(N^{-2/3})$ is the same as the one of the fluctuations of the
largest eigenvalue of the complex Wishart model is crucial.

We detail the proof in the case when $\rho > \sqrt c$ and, as above,
consider the positive measures $\mathbb Q_1$. We need to prove that:
\begin{equation}\label{eq:lower-precise}
\liminf_{N,K\rightarrow\infty} \frac{1}{N} \log  \mathbb Q_1\left(T_N < (1+\sqrt{c_N})^2 + \frac{\eta}{N^{2/3}}\right)
 \geq - G_\rho(\lambda^+),\quad \eta\in \mathbb R,
\end{equation}
the other bound being a direct consequence of the LDP. As previously,
we will carefully localize the various quantities of interest.
Denote by $g_N(\eta) = (1+\sqrt{c_N})^2 + \eta N^{-2/3}$ for $\eta \in
\mathbb R$ and by $h_N(r) = 1- r N^{-2/3}$ for $r>0$. Notice also that
$\lambda_1 \leq g_N(\eta)h_N(r)$ together with $\frac{1}{K-1}\sum_{j=2}^{K}\lambda_j > h_N(r)$ imply that $T_N < g_N(\eta)$.
We shall also consider the further constraints:
$$
g_N(\eta-1) h_N(r) \leq \lambda_1\quad \textrm{and}\quad \lambda_2 <
g_N(\eta-2)h_N(r)
$$
which enable us to properly separate $\lambda_1$ from the support of  $\hat \pi_{K,\boldsymbol{\lambda}}$.
Now, with the localisation indicated above, we have
for $N$ large enough,
\begin{multline*}
 \mathbb Q_1\left(T_N < g_N(\eta)\right)
\ \geq\ \mathbb Q_1\bigg( g_N(\eta-1) h_N(r)
\leq \lambda_1 \leq g_N(\eta)h_N(r), \\ \frac{1}{K-1}\sum_{j=2}^{K}\lambda_j > h_N(r),\ \lambda_2 < g_N(\eta-2)h_N(r),\
\hat \pi_{K,\boldsymbol{\lambda}} \in B( \Pmp,N^{-1/4})\bigg).
\end{multline*}
As previously, we consider the variables
$y_j = \frac{N}{N-1}x_j$ for $2\leq j\leq K$ and obtain, with the help of Lemma \ref{IZlemma}:
\begin{multline*}
\mathbb Q_1\left(T_N < g_N(\eta)\right)
\geq \int_{ g_N(\eta-1) h_N(r)}^{g_N(\eta) h_N(r)} dx_1 \int_{ \mathcal F}
e^{-N \Phi(x_1,c_N, \hat \pi_{K,\bf y})} e^{Nc(J_{\rho}(x_1)- \delta)} dp_{K-1,N-1}^{0}(y_{2:K})
\end{multline*}
with
\begin{multline*}
\mathcal F := \left\{(y_2,\cdots,y_K)\in \left[0,
    \frac{N\,g_N(\eta-2)\,h_N(r)}{N-1}\right]^{K-1},\right.\\
\left.   \frac{1}{K-1}\sum_{j=2}^{K}y_j > \frac{N}{N-1}h_N(r), \hat
  \pi_{K,\bf y} \in B( \Pmp,N^{-1/4}) \right\}.
\end{multline*}
Therefore:
$$
\mathbb Q_1\left(T_N < g_N(\eta)\right) \ \geq\ h_N(r)\ (g_N(\eta) - g_N(\eta-1))\ e^{N( G_\rho(\lambda^+) -
  2\delta)}\mathbb P_0\left( (\lambda_2,\cdots,\lambda_K)\in \mathcal
  F\right)
$$
(recall that $G_\rho(x) = \Phi(x,c, \Pmp) + c J_\rho\left(
  x\right)$). Now, as $h_N(r)\ (g_N(\eta) - g_N(\eta-1))=(1-rN^{-2/3})N^{-2/3}$, its contribution
vanishes at the LD scale:
$$
\lim_{N\rightarrow\infty} \frac 1N \log \left( h_N(r)\ (g_N(\eta) - g_N(\eta-1))\right) =0\ .
$$ 
It remains to check that $\mathbb
P_0\left((\lambda_2,\cdots,\lambda_K)\in \mathcal F\right)$ is bounded
below uniformly in $N$. This will yield the convergence of
$\frac{1}{N} \log \mathbb P_0\left(
  (\lambda_2,\cdots,\lambda_K)\in\mathcal F\right)$ towards zero,
hence \eqref{eq:lower-precise}. Consider:
\begin{multline*}
\mathbb P_0\left( (\lambda_2,\cdots,\lambda_K)\in \mathcal F^c\right) \leq
\mathbb P_0\left( \hat \pi_{K,{\boldsymbol \lambda}} \notin B( \Pmp,N^{-1/4})\right)\\
+\mathbb P_0\left( \frac{1}{K-1}\sum_{j=2}^{K}\lambda_j < \frac{N}{N-1}h_N(r)\right)
+ \mathbb P_0\left( \lambda_2 > \frac{N}{N-1}g_N(\eta-2)h_N(r)\right)\ .
\end{multline*}
We have already used the fact that the first term goes to zero when
$N$ grows to infinity.  Recall that the fluctuations of
$\frac{1}{K-1}\sum_{j=2}^{K}\lambda_j $ are of order $\frac{1}{N}$,
therefore the second term also goes to zero as we consider deviations
of order $N^{-2/3}$.  Now, $N^{2/3}(\lambda_2 -
(1+\sqrt{c_N})^2)$ converges in distribution to the Tracy-Widom law,
therefore the last term converges to $F_{\rm TW}\left( \eta-2 +
  r(1+\sqrt{c})^2\right) <1.$ This concludes the proof.

\section{Sketch of proof for Lemma~\ref{lem:ldpU}: Large deviations for $U_N$}
\label{app:ldp-U}
As stated in Remark \ref{rem:lem2}, we shall first study the LDP for
the joint quantity $(\lambda_1,\lambda_K)$. The purpose here is to outline the
following convergence:
$$
\frac 1N \log \mathbb{P}\left(\lambda_1\in A,\lambda_K\in B\right)
\xrightarrow[N,K\to \infty]{} - \inf_{x\in A} I^+_{\rho}(x) - \inf_{y\in B} I^-(x)\ , 
$$
which is an illustrative way, although informal\footnote{All the statements,
computations and approximations below can be made precise as in the
proof of Lemma \ref{lem:ldp}.}, to state the LDP for
$(\lambda_1,\lambda_K)$ (see \eqref{eq:ldp-informal}).

Consider the quantity $\mathbb{P}\left(\lambda_1\in
  (\alpha_1,\beta_1),\lambda_K \in (\alpha_K,\beta_K)\right)$. As we
are interested in the deviations of $\lambda_1$ and $\lambda_K$, the
interesting scenario is $\lambda^+\notin (\alpha_1,\beta_1)$ and
$\lambda^-\notin (\alpha_K,\beta_K)$ (recall that $\lambda^\pm$ are
the edgepoints of the support of Mar$\check{\textrm{c}}$enko-Pastur distribution). More
precisely, the interesting case is when the deviations of the extreme
eigenvalue occur outside of the bulk: $\alpha_1 > \lambda^+$ and
$\beta_K< \lambda^-$; such deviations happen at the rate $e^{-N \times
  const.}$. The case where the deviations would occur within the bulk
is unlikely to happen because it would enforce the whole eigenvalues
to deviate from the limiting support of Mar$\check{\textrm{c}}$enko-Pastur distribution,
which happens at the rate $e^{-N^2 \times const.}$. Denote by
$A=(\alpha_1,\beta_1)$ and $B=(\alpha_K,\beta_K)$. 
\begin{eqnarray*}
  \lefteqn{  \mathbb{P}\left(\lambda_1\in A,\lambda_K\in B\right)}\\
  &=&\frac 1{Z_{K,N}^1} \int_{A\times \mathbb{R}^{(K-2)}\times B}
  1_{(\lambda_1\ge \cdots\ge \lambda_K \ge 0)} 
  \prod_{1\le i< j\le K} (x_i - x_j)^2 \\
  &&\qquad \qquad \qquad \times \prod_{j=1}^K x_j^{N-K} e^{-N x_j} I_K\left( \frac NK B_K,X_K\right)
  d\,x_{1:K} \\
  &=& \int_A d\, x_1\  e^{2\sum_{j=2}^{K-1} \log(x_1 - x_j)} e^{(N-K)\log x_1 -N x_1} I_K\left( \frac NK B_K,X_K\right)\\
  && \quad \times \int_B d\, x_K\  e^{2\sum_{i=2}^{K-1} \log(x_i - x_K)} e^{(N-K)\log x_K -N x_K} e^{2\log(x_1 - x_K)}\\
  &&\quad \quad \times \frac{Z_{K-2,N-2}^0}{Z_{K,N}^1} \int_{x_1\ge x_2\ge \cdots \ge x_K} 
  \prod_{j=2}^{K-1} e^{-2 x_j}  
  \prod_{j=2}^{K-1} \frac{x_j^{N-K} e^{-(N-2)x_j}}{Z_{K-2,N-2}^0} \prod_{2\le i<j\le K-1} (x_i - x_j)^2 d\, x_{2:K-1}
\end{eqnarray*}
We shall now perform the following approximations:
\begin{eqnarray*}
\sum_{j=2}^{K-1} \log(x_1 -x_j) &\approx& (K-2)\int \log(x_1 -x)\Pmp(\,dx)\ = (K-2) {\bf F}^+(x_1)\ ,\\
\sum_{j=2}^{K-1} \log(x_j -x_K) &\approx& (K-2)\int \log(x -x_K)\Pmp(\,dx)\ = (K-2) {\bf F}^-(x_K)\ ,\\
\sum_{j=2}^{K-1} x_j &\approx& (K-2)\int x \Pmp(\,dx)\ = (K-2) \ ,\\
I_K\left( \frac NK B_K,X_K\right) &\approx& e^{N c J_{\rho}(x_1)}\ .
\end{eqnarray*}
The three first approximations follow from the fact that $\frac 1{K-2}
\sum_2^{K-1} \delta_{x_i} \approx \Pmp$, the last one from Lemma
\ref{IZlemma}.  Plugging these approximations into the expression of
$\mathbb{P}\left(\lambda_1\in A,\lambda_K\in B\right)$ yields:
\begin{eqnarray*}
\lefteqn{  \mathbb{P}\left(\lambda_1\in A,\lambda_K\in B\right)}\\
  &\approx& \int_A d\, x_1\  e^{2(K-2) {\bf F}^+(x_1)}  e^{(N-K)\log x_1 -N x_1} e^{N cJ_{\rho}(x_1)}\\
  && \quad \times \int_B d\, x_K\  e^{2(K-2) {\bf F}^- (x_K)} e^{(N-K)\log x_K -N x_K} e^{2\log(x_1 - x_K)}\\
  &&\quad \quad \times \frac{Z_{K-2,N-2}^0}{Z_{K,N}^1} e^{-2(K-2)} \int_{x_1\ge x_2\ge \cdots \ge x_K} 
\prod_{j=2}^{K-1} \frac{x_j^{N-K} e^{-(N-2)x_j}}{Z_{K-2,N-2}^0} \prod_{2\le i<j\le K-1} (x_i - x_j)^2 d\, x_{2:K-1}\ .
\end{eqnarray*}
As $x_1\ge \alpha_1\ge \lambda^+$ and $x_K\le \beta_K \le \lambda^-$,
the last integral goes to one as $K,N\to \infty$ and:
\begin{eqnarray*}
\lefteqn{  \mathbb{P}\left(\lambda_1\in A,\lambda_K\in B\right)}\\
  &\approx& \int_A d\, x_1\ 
e^{
-N\left( \frac{2(K-2)}N {\bf F}^+(x_1) - \left( 1 -\frac KN\right) \log x_1 + x_1 -cJ_{\rho}(x_1) \right)
}\\
&&\times \int_B d\, x_K\ 
e^{
-N\left( \frac{2(K-2)}N {\bf F}^-(x_K) - \left( 1 -\frac KN\right) \log x_K + x_K  + \frac{2\log(x_1 - x_K)}N\right)
} \\
  && \quad \times \frac{Z_{K-2,N-2}^0}{Z_{K,N}^1} e^{-2(K-2)} \ .
\end{eqnarray*}
Recall that we are interested in the limit $N^{-1} \log
\mathbb{P}\left(\lambda_1\in A,\lambda_K\in B\right)$. The last term
will account for a constant $\Upsilon$ (see for instance
\eqref{eq:conv-constante}):
$$
\frac 1n \log \left( \frac{Z_{K-2,N-2}^0}{Z_{K,N}^1} e^{-2(K-2)}
\right)\quad \xrightarrow[N,K\to \infty]{} \quad \Upsilon\ .
$$
The term $\frac{2\log(x_1 - x_K)}N$ within the exponential in the
integral accounts for the interraction between $\lambda_1$ and
$\lambda_K$ and its contribution vanishes at the desired rate. In
order to evaluate the two remaining integrals, one has to rely on
Laplace's method (see for instance \cite{Die71})  to express the leading term of the
integrals (replacing $KN^{-1}$ by $c$ below):
\begin{eqnarray*}
  \int_A d\, x_1\ e^{
    -N\left( 2c {\bf F}^+(x_1) - \left( 1 -c\right) \log x_1 + x_1 -cJ_{\rho}(x_1) \right)
  }&\approx& e^{-N\inf_{x\in A} \left( 2c {\bf F}^+(x) - \left( 1 -c\right) \log x + x -cJ_{\rho}(x) \right) }\ ,\\
\int_B d\, x_K\ 
e^{
-N\left( 2c {\bf F}^-(x_K) - \left( 1 -c \right) \log x_K + x_K -cJ_{\rho}(x_K) \right)}&\approx& 
e^{-N\inf_{y\in B} \left( 2c {\bf F}^-(y) - \left( 1 -c\right) \log y + y  \right) }\ .
\end{eqnarray*}
Finally, we get the desired limit:
$$
\frac 1N \log \mathbb{P}\left\{\lambda_1\in A,\lambda_K\in B\right\}
\xrightarrow[N,K\to \infty]{} - \inf_{x\in A} \Phi^+(x) - \inf_{y\in B} \Phi^-(y) + \Upsilon\ ,
$$
where
\begin{eqnarray*}
\Phi^+(x)&=& 2c {\bf F}^+(x) - \left( 1 -c\right) \log x + x -cJ_{\rho}(x) \ , \\
\Phi^-(y)&=& 2c {\bf F}^-(y) - \left( 1 -c\right) \log y + y \ .
\end{eqnarray*}
It remains to replace $J_{\rho}$ by its expression \eqref{eq:jrho} and
to spread the constant $\Upsilon$ over $\Phi^+$ and $\Phi^-$, which
are not a priori rate functions (recall that a rate function is
nonnegative). If $\lambda^- \in B$, then the event $\{\lambda_K \in
B\}$ is ``typical'' and no deviation occurs, otherwise stated, the
rate function $I^-$ should satisfy $I^-(\lambda^-)=0$. Similarly,
$I_0^+(\lambda^+)=0$ under $H_0$ and $I_{\rho}^+(\lspiked)=0$ under
$H_1$. Necessarily, $\Upsilon$ should write $\Upsilon=
\Phi(\lambda^-)+\Phi(\lambda^+)$ under $H_0$ (resp. $\Upsilon=
\Phi(\lambda^-)+\Phi(\lspiked)$ under $H_1$) and the rate functions
should be given by: $I^-= \Phi^- - \Phi(\lambda^-)$, $I_0^+= \Phi^+ -
\Phi(\lambda^+)$ under $H_0$ (resp.  $I_\rho^+=\Phi^+ -\Phi(\lspiked)$
under $H_1$), which are the desired results.

We have proved (informally) that the LDP holds true for
$(\lambda_1,\Lambda_K)$ with rate function $I_{0/\rho}^+(x)+I^-(y)$. The
contraction principle \cite[Chap. 4]{DemZei98} immediatly yields the
LDP for the ratio $\frac{\lambda_1}{\lambda_K}$ with rate function:
\begin{equation}\label{eq:rf-ratio}
\Gamma_{0/\rho}(t) = \inf_{(x,y), \frac xy =t}\{ I^+_{0/\rho}(x) + I^-(y)\} \ ,
\end{equation}
which is the desired result. We provide here intuitive arguments to
understand this fact. 

For this, interpret the value of the rate function $I_\rho^+(x)$ as
the cost associated to a deviation of $\lambda_1$ (under $H_1$) around
$x$: $\mathbb{P}\{\lambda_1 \in (x,x+dx)\} \approx
e^{-NI^+_{\rho}(x)}$. If a deviation occurs for the ratio
$\frac{\lambda_1}{\lambda_K}$, say $ \frac{\lambda_1}{\lambda_K} \in
(t , t+dt)$ where $t>\frac{\lspiked}{\lambda^-}$ (which is the typical
behaviour of $U_N$ under $H_1$), then necessarily $\lambda_1$ must
deviate around some value $ty$, so does $\lambda_K$ around some value
$y$, so that the ratio is around $t$. In terms of rate functions, the
cost of the joint deviation $(\lambda_1\approx ty, \lambda_K \approx
y)$ is $I^+_\rho(ty) + I^-(y)$. The true cost associated to the
deviation of the ratio will be the minimum cost among all these
possible joint deviations of $\lambda_1$ and $\lambda_K$, hence the
rate function \eqref{eq:rf-ratio}.

\section{Closed-form expressions for functions $\f$, $\mathbf{F}^+$ and $\mathbf{F}^-$}
\label{appendix:representation}
Consider the Stieltjes transform $\f$ of Mar$\check{\textrm{c}}$enko-Pastur distribution:
$$
\f(z)=\int \frac{\Pmp(d\lambda)}{\lambda-z}\ .
$$
We gather without proofs a few facts related to $\f$,
which are part of the folklore.
\begin{lemma}[Representation of $\f$]\label{folklore-f} The following hold true:
\begin{enumerate}
\item Function ${\bf f}$ is analytic in $\mathbb{C} - [\lambda^-,\lambda^+]$.
\item If $z\in \mathbb{C} - [\lambda^-,\lambda^+]$ with $\Re(z)\geq \frac{\lambda^+ + \lambda^-}2$, then
$$
{\bf f}(z)= \frac{(1-z-c) + \sqrt{(1-z-c)^2 - 4cz}}{2cz}\ ,
$$
where $\sqrt{z}$ stands for the principal branch of the square-root.
\item If $z\in \mathbb{C} - [\lambda^-,\lambda^+]$ with $\Re(z) < \frac{\lambda^+ + \lambda^-}2$, then
$$
{\bf f}(z)= \frac{(1-z-c) - \sqrt{(1-z-c)^2 - 4cz}}{2cz}\ ,
$$
where $-\sqrt{z}$ stands for the branch of the square-root whose image is $\{z\in \mathbb{C},\ \Re(z)\leq 0\}.$
\item As a consequence, the following hold true:
\begin{eqnarray}
\quad {\bf f}(x) &=& \frac{(1-x-c) + \sqrt{(1-x-c)^2 - 4cx}}{2cx}\quad \textrm{if} \quad x\geq\lambda^+\ ,\label{eq:fplus}\\
\quad {\bf f}(x) &=& \frac{(1-x-c) - \sqrt{(1-x-c)^2 - 4cx}}{2cx}\quad \textrm{if} \quad 0\leq x\leq \lambda^-\ .\label{eq:fmoins}
\end{eqnarray}
\item Consider the following function
$
\tildef(z)= c \f(z) -\frac{1-c}{z}
$.
Functions $\f$ and $\tildef$ satisfy the following system of equations:
\begin{equation}\label{ftildef}
\left\{
\begin{array}{lcl}
\f(z)&=& -\frac{1}{z(1+\tildef(z))}\\
\tildef(z)&=& -\frac{1}{z(1+c\f(z))}
\end{array}
\right. \ ,
\end{equation}

\end{enumerate}

\end{lemma}

Recall the definition \eqref{eq:Fp} and \eqref{eq:Fm} of function $\mathbf{F}^+$
and $\mathbf{F}^-$. In the following lemma, we provide closed-form formulas of interest.
\begin{lemma} \label{lem:F} The following identities hold true:
\begin{enumerate}
\item Let $x\geq\lambda^+,$ then
$$
  \mathbf{F}^+(x) = \log(x) +\frac 1c \log(1+c \mathbf{f}(x)) +\log(1+\tilde{\mathbf{f}}(x)) +x\mathbf{f}(x)
\tilde{\mathbf{f}}(x)\ .
$$
\item Let $0\leq x\leq\lambda^-,$ then
$$
\mathbf{F}^-(x) =
\log(x) +\frac 1c \log(1+c \mathbf{f}(x)) +\log(-(1+\tilde{\mathbf{f}}(x)))+x\mathbf{f}(x)
\tilde{\mathbf{f}}(x)\ .
$$
\end{enumerate}
\end{lemma}

\begin{proof} Consider the case where $x\geq \lambda^+$.
First write
$$
\log(x-y) = \log(x) +\int_x^{\infty} \left( \frac 1u + \frac 1{y-u}\right) \,du \ .
$$
Integrating with respect with $\Pmp$ and applying Funini's theorem yields:
$$
\int \log(x-y) \Pmp(\,dy) = \log(x) + \int _x^\infty \left( \frac 1u + \f(u) \right) \,du
$$
in the case where $x>\lambda^+$. Recall that $\f$ and $\tildef$ are
holomorphic functions over $\mathbb{C} - (\{0\}\cup [\lambda^-,\lambda^+])$
and satisfy system \eqref{ftildef} (notice in particular that $1+c\f$
and $1+\tildef $ never vanish). Using the first equation of \eqref{ftildef} implies that:
\begin{equation}\label{eq-integ-intermediaire}
\int \log(x-y) \Pmp(\,dy) = \log(x) - \int _x^\infty \f(u) \tildef(u) \,du\ .
\end{equation}
Consider $\Gamma(u,\f, \tildef)= \frac 1c \log (1+c \f) +\log
(1+\tildef) +u\f \tildef$. By a direct computation of the derivative,
we get:
\begin{eqnarray*}
\frac{d}{du} \Gamma(u,\f(u),\tildef(u)) &=& \f'\left( \frac 1{1+c\f} +u\tildef\right)
+ \tildef'\left( \frac 1{1+\tildef} +u\f\right) +\f \tildef\\
&=& \f(u) \tildef(u)\ .
\end{eqnarray*}
Hence
\begin{eqnarray*}
\int_x^\infty \f(u) \tildef(u) \,du &=&\left[ \frac 1c \log (1+c \f) +\log
(1+\tildef) +u\f \tildef \right]_x^\infty\\
&=&  - \left( \frac 1c \log (1+c \f(x) ) +\log
(1+\tildef(x) ) +x\f(x) \tildef(x) \right).
\end{eqnarray*}
It remains to plug this identity into \eqref{eq-integ-intermediaire} to conclude.
The representation of $\mathbf{F}^-$ can be established similarly.

\end{proof}

\bibliographystyle{unsrt}
\bibliography{eigen,math}
\end{document}